\documentclass[]{amsart}
\usepackage{amsmath}
\usepackage{amssymb}
\usepackage{anysize}

\marginsize{2.8cm}{2.8cm}{2.8cm}{2.8cm}

\theoremstyle{plain}
\newtheorem{theor}{Theorem}[section]
{}
\newtheorem{lemma}{Lemma}[section]
\newtheorem{propo}{Proposition}[section]
\newtheorem{corol}{Corollary}[section]

\theoremstyle{definition}
\newtheorem{defin}{Definition}[section]

\newtheorem{notat}{Notation}[section]
\newtheorem{quest}{Question}[section]

\newtheorem*{case1}{Case One}
\newtheorem*{case2}{Case Two}

\newtheorem{claim}{Claim}

\theoremstyle{remark}

\numberwithin{equation}{section}

\DeclareMathOperator{\dom}{domain}

\newcommand{\forcingname}[1]{\buildrel \circ\over #1}
\newcommand{\AND}{\text{ and }}
\newcommand{\IF}{\text{ if }}
\newcommand{\OR}{\text{ or }} 
\newcommand{\card}[1]{\lvert #1 \rvert}

\newcommand{\CH}{$2^{\aleph_0} = \aleph_1$\ }
\newcommand{\forces}[2]{\Vdash_{#1} \mbox{``} #2 \mbox{''}}

\newcommand{\Sacks}{{\mathbb S}}
\newcommand{\Laver}{{\mathbb L}}
\newcommand{\Reals}{{\mathbb R}}

\newcommand{\Poset}{{\mathbb P}}

\newcommand{\Integers}{{\mathbb Z}}
\newcommand{\Naturals}{{\mathbb N}}

\newcommand{\Fin}{{\mathbb F}}

\newcommand{\presup}[2]{\, ^{#1} \! #2} 
\newcommand{\wpresup}[1]{\presup{\stackrel{\omega}{\smile}}{#1}} 
\newcommand{\fomom}{\presup{\omega}{\omega}}
\newcommand{\wfomom}{\wpresup{\omega}}
\newcommand{\pomega}{{\mathcal P}(\omega)}
\newcommand{\pomegaf}{{\mathcal P}(\omega)/[\omega]^{<\aleph_0}}
\title[Maximal Abelian Subgroups]{Possible Cardinalities of Maximal
  Abelian Subgroups of Quotients of Permutation Groups of the Integers}
\author[S. Shelah]{Saharon Shelah}
\address{Department of Mathematics, Rutgers University, Hill Center,
 Piscataway, 
 New Jersey, U.S.A. 08854-8019}
\curraddr{Institute of Mathematics\\Hebrew University\\
Givat Ram, Jerusalem 91904, Israel}
\email{shelah@math.rutgers.edu}
\author[J. Stepr\={a}ns]{Juris Stepr\={a}ns}
\address{Department of Mathematics, York University\\
4700 Keele Street,
North York, Ontario\\ Canada \ \ \  M3J 1P3}
\curraddr{}
\email{steprans@yorku.ca}
\thanks{The research of the first author was supported by The Israel Science
Foundation founded by the Israel Academy of Sciences and Humanities, and
by NSF grant No. NSF-DMS97-04477.
Research of the second author for this paper was partially supported by NSERC of
Canada. This is number 786 on the first author's personal
list of publications. }
\begin{document}
\maketitle
\bibliographystyle{plain}
\section{Introduction and Definitions}
The maximality of Abelian subgroups play a role in various parts of
group theory. For example,
Mycielski \cite{MR20:6479,MR20:1725}
  has extended a classical result of Lie groups and  shown that
 a maximal   Abelian subgroup of a compact connected group is
 connected and, furthermore, all
   the maximal Abelian subgroups are conjugate.
For finite symmetric groups the question of the size of maximal Abelian
subgroups has been examined by
Burns and Goldsmith in \cite{MR89i:20013} and
 Winkler in \cite{MR94m:20013}.
It will be shown
 in Corollary~\ref{c:ps} that there is not much interest in
 generalizing this study to infinite symmetric groups;
the cardinality
 of any
 maximal Abelian subgroup of the symmetric group of the integers 
is $2^{\aleph_0}$. 
The purpose of this paper is to examine the size of maximal Abelian
subgroups for a class of groups closely related to the
the symmetric group of the integers; these arise by taking an ideal on
the integers, considering
the subgroup of all permutations which respect the ideal and then
taking the quotient by the normal subgroup of permutations which fix
all integers except a set in the ideal. It will be shown that the
maximal size of Abelian subgroups in such groups is sensitive to 
the nature of the ideal as well as various set theoretic hypotheses. 

The reader familiar with applications of the  Axiom of Choice may not
be surprised 
by the assertion just made since, one can imagine constructing ideals
on the integers by transfinite induction such that the quotient group just
described exhibits various desired properties. Consequently, it is of
interest to  restrict attention to
 only those ideals which do not require the Axiom
of Choice for their definition.
All of the ideals considered will here will have
 simple definitions --- indeed, they will all be
Borel subsets of $\pomega$ with the usual topology --- and, in fact,
the first three sections will focus on the ideal of finite sets. 
It should be
mentioned that there is large body of work examining 
the analogous quotients of the Boolean algebra $\pomega$ modulo an 
  analytic ideal --- the monograph \cite{MR2001c:03076} by Farah is a
  good reference for this subject. However, the analogy is far from
  perfect since, for example,
whereas the Boolean algebra $\pomegaf$ can consistently have $2^{2^{\aleph_0}}$
automorphisms \cite{MR18:324d} it is shown in  \cite{MR99a:20002} that 
the quotient of the full symmetric group of the integers 
modulo the subgroup of finite permutations has 
only countably many outer automorphisms. Nevertheless, it may be
possible to employ methods similar to those of  \cite{MR2001c:03076}
in order to distinguish between different quotient algebras up to
isomorphism. This has been done for elementary equivalence
in \cite{MR2000f:03114,MR54:4972} for quotients of the full symmetric group
on $\kappa$ by the normal subgroups fixing all but $\lambda$ elements.
However since the full symmetric group of the integers has only two proper
normal subgroups \cite{sc.ul.1933}
quotients of certain naturally arising subgroups will
be considered instead. One of the goals of this study is to use the
cardinal invariant associated with maximal Abelian subgroups as a
tool to distinguish between isomorphism types of such groups.

In order to state the main results
  precisely some notation is needed.
\begin{defin}
  If $G$ is a group then define the Abelian subgroup spectrum of $G$ to be the
   set of all $\kappa$ such that there is a maximal Abelian subgroup of
   $G$ of 
  size $\kappa$. Define $A(G)$ to be least uncountable cardinal in the
   Abelian subgroup spectrum of $G$.  
\end{defin}
\begin{notat}
  Through this paper the symbol $\Sacks$ will be used to denote the
  symmetric 
  group on  $\Naturals$. 
  For $ \pi\in \Sacks$ let $\text{supp}(\pi)$ denote
  the support of $\pi$ which is defined to be $\{x \in \dom(\pi): \pi(x)\neq
  x\}$. 
If ${\mathcal I}$ is an ideal\footnote{An ideal is a collection of
  subsets of the integers closed under finite unions and subsets.}
 on $\Naturals$
  then $\Sacks({\mathcal I})\subseteq \Sacks$ will denote the subgroup
  of all permutations preserving ${\mathcal I}$; in other words,
a permutation $\pi$ belongs to $ \Sacks({\mathcal I})$ provided 
$\pi(A) \in {\mathcal I}$ if and only if $A\in {\mathcal I}$.
   On the other hand,  $\Fin({\mathcal I})$ will be used to denote the normal 
  subgroup of $\Sacks({\mathcal I})$
 consisting of all permutations $\pi \in\Sacks({\mathcal I})$
such that $\text{supp}(\pi)\in {\mathcal I}$. The abbreviation
$\Fin = \Fin([\Naturals]^{<\aleph_0}$ will also be used.
\end{notat}

The focus of this paper will be on 
examining 
computing $A(\Sacks({\mathcal
  I})/\Fin({\mathcal I}))$ for various simply defined  ideals.
 This cardinal will be denoted by $A({\mathcal I})$.

\begin{notat}
  Given a pair of permutations $\{\pi,\pi'\}\in [\Sacks]^2$ define
$\text{NC}(\pi,\pi') = \{n\in \Naturals : 
\pi(\pi'(n)) \neq \pi'(\pi(n))\}$. A pair 
of permutations $\{\pi,\pi'\}\in [\Sacks]^2$ will be said to almost commute modulo an ideal $\mathcal I$ if 
$\text{NC}(\pi,\pi')\in \mathcal I$ and they will be said 
  to almost commute if
$\text{NC}(\pi,\pi')$ is finite. 
\end{notat}
\begin{notat}
  Given a  permutation $\pi$ and a $X\subseteq \Naturals$
define the orbit
of $X$ under $\pi$
by $\text{orb}_{\pi}(X) = \{\pi^i(x)\}_{i\in \Integers,x\in X}$.
The abbreviation $\text{orb}_{\pi}(n) = \text{orb}_{\pi}(\{n\})$ will
be used when no confusion is possible.  
\end{notat}

\begin{notat}
  Given a permutation $\pi\in \Sacks$ define $\equiv_{\pi}$  be the
  equivalence relation on $\Naturals$ whose equivalence classes are the orbits
  of $\pi$. 
  Given a set of permutations ${\mathcal S}\subseteq \Sacks$ define
  $\equiv_{\mathcal S}$ to be the transitive closure of the set of equivalence
  relations $\{\equiv_{\pi}\}_{\pi\in {\mathcal S}}$. 
Let $\Omega_{\mathcal S}$ denote the set of equivalence classes of
$\equiv_{\mathcal S}$ and, 
for any set $X$ define $\Omega_{\mathcal S}(X)$ to be
 the smallest set containing $X$ and closed
  under equivalence classes of $\equiv_{\mathcal S}$.
\end{notat}

\begin{notat}\label{d:324}
Given two finite subsets $A$ and $B$ of $\Naturals$ define
$\Delta_{A,B}:A \to B$ to be the unique order preserving mapping
between them and let $\Delta_A = \Delta_{A, |A|}$.
If ${\mathcal S} \subseteq \Sacks$
and $A$ and $B$ are two
   equivalence classes of $\equiv_{\mathcal S}$ then define $A$ to be
   ${\mathcal S}$-{\em
   isomorphic} 
   to $B$ if there is a bijection $\psi:A\to B$ such that $\pi(\psi(a)) =
   \psi(\pi(a))$ for each $\pi \in\mathcal S$ and $a\in A$.  
\end{notat}
The set theoretic notation used throughout will follow the contemporary
standard. In particular, $[X]^k$ will denote the family of subsets of
$X$ of cardinality $k$ and $[X]^{<k}$ will denote the family of subsets of
$X$ of cardinality less than $k$. Cardinal invariants of the continuum
are closely linked to investigation of $A({\mathcal I})$. The
following recall the definitions of some well known invariants.
\begin{defin}\label{d:d12}
Given an ideal ${\mathcal I}\subseteq \pomega$
let $\pomega/{\mathcal I}$ be the quotient Boolean algebra
 and denote the least cardinal of a maximal, uncountable, pairwise disjoint
 family\footnote{See \cite{step.49} for a more detailed discussion
  of this invariant.} in $\pomega/{\mathcal I}$ is denoted by ${\mathfrak
 a}({\mathcal I})$.
In the special case ${\mathcal I} = [\Naturals]^{<\aleph_0}$
${\mathfrak
 a}({\mathcal I})$ is denoted by ${\mathfrak
 a}$
The least cardinal of an ideal 
   ${\mathcal B} \subseteq
  \pomegaf$ such that there is no $C \in \pomegaf$ disjoint from all
  members of ${\mathcal B}$ other than the equivalence class of the
  finite sets is denoted by $\mathfrak p$.
\end{defin}

In Section~2 it is shown that ${\mathfrak a}$  is an upper bound
for $A([\Naturals]^{<\aleph_0})$ while in Section~3 it is shown that
 ${\mathfrak p}$ serves as a lower bound
for $A([\Naturals]^{<\aleph_0})$. Sections~4 and~5 deal with
consistency results. In Section~4 it is shown that ${\mathfrak a}$ is
not the best possible upper bound for $A([\Naturals]^{<\aleph_0})$
since in the iterated Laver model $A([\Naturals]^{<\aleph_0})$ is
strictly less than ${\mathfrak a}$. 
Sections~5 and ~6 deal with quotients using ideals other than the
ideal of finite sets. It is shown in Section~5 that adding $\aleph_1$
Cohen reals to a model where 
$2^{\aleph_0} > \aleph_1$ yields a model where 
$A({\mathcal I}_{1/x}) = \aleph_1 < 2^{\aleph_0}$ and ${\mathcal
  I}_{1/x}$ is the ideal of sets whose reciprocals form a series with
finite sum.
Section~6 deals with ideals similar to the density ideal. It
is shown that $A({\mathcal I}) = 2^{\aleph_0}$
for many of these ideals $\mathcal I$. No extra set theoretic axioms
are used here. The final section contains some open questions.

\section{An upper bound}
\begin{propo}\label{p:1}
$A([\Naturals]^{<\aleph_0})
\leq {\mathfrak a}$.
\end{propo}
\begin{proof}
 Let $\mathcal A$ be a maximal almost disjoint family of subsets of
  $\Naturals$  of size $\mathfrak a$ and let $F({\mathcal A})$ be the
  free Abelian group generated by $\mathcal A$ under coordinate wise addition; in other words,  if $f
  \in F({\mathcal A})$ then  $f:{\mathcal A} \to \Integers$ and $f$
  has finite support.
For $a \in {\mathcal A}$ define $\pi_a: a \to a$ by
$\pi_a(i) = 
 \min(\{j \in a : j > i\})$ and, for $j \in \Integers$, let $\pi_a^j$
  denoted the $j$-fold composition of $\pi_a$ noting that both the
  domain and range of $\pi^j_a$ are co-finite subsets of $a$.. 
If $f \in F({\mathcal A})$ then  let $\Phi(f)$ be the set of all
permutations $\pi$ such that there is a finite set $F\subseteq \Naturals$
such that
$$\pi(j) = 
\begin{cases}
\pi_a^{f(a)}(j) & \IF j \in a\setminus F \AND f(a)\neq 0\\
j & \IF j \notin F\AND (\forall a\in {\mathcal A})j\notin a \OR f(a)
= 0  
\end{cases}$$
leaving $\Phi(f)$ undefined if there are no such permutations.
Note that if $\pi \in \Phi(f)$ and $\sigma \in \Phi(g)$ then
$\pi\circ\theta\in \Phi(f+g)$.
Since $\Phi$ is easily seen to be one-to-one, it is an 
isomorphism between a subgroup of $F({\mathcal A})$ and the subgroup
$\Phi(F({\mathcal
    A}))$ of $\Sacks/\Fin$. 

In fact, $\Phi(f)$ is defined precisely when $\sum_{a\in {\mathcal A}}f(a)
= 0$. To see this, let $f \in F({\mathcal A})$ and suppose that the
support of $f$ is $B$ and $\sum_{b\in B}f(b) = 0$. Let $F \subseteq
\Naturals$ be a finite set such that $b\cap b'\subseteq F$ for any two
$b$ and $b'$ in $B$ and such that $F\cap b$ is an initial segment of
$b$ for each $b\in B$.  Let $B^+ = \{b\in B : f(b) > 0\}$ and $B^- =
\{b\in B : f(b) < 0\}$.
If $b \in B^+$ let $b^*$ be the first $f(b)$ elements of
$b\setminus F$ and if $b\in B^-$ let $b^*$ be the first $-f(b)$ elements of
$b\setminus F$. Let $\theta: \bigcup_{b\in B^-}b^* \to \bigcup_{b\in
  B^+}b^*$ be any bijection and define $\pi$ as follows:
$$
\pi(j) =
\begin{cases}
  j & \IF j \notin \bigcup_{b\in B} b\setminus F\\
\pi_b^{f(b)} & \IF j \in b \in B^+\\
\pi_b^{f(b)} & \IF j \in b\setminus b^*\AND b \in B^-\\
\theta(j) & \IF j \in \bigcup_{b\in B^-}b^*
\end{cases}$$
and observe that $\pi$ is a bijection.

To see that $\Phi(F({\mathcal
    A}))$ is maximal Abelian let $[\pi]_{\Fin} \in \Sacks/\Fin
\setminus\Phi(F({\mathcal
    A})) $. 
Before continuing some notation will be introduced. Given two distinct
    elements $a$ and $a'$ of ${\mathcal
    A}$ define
$f_{a,a'} \in F({\mathcal
    A})$ be such that $\text{supp}(f_{a,a'}) = \{a, a'\}$ and $f_{a,a'}(a) = 1 =
  -f_{a,a'}(a')$. 
 Choose $\pi_{a,a'}\in \Phi(f_{a,a'})$.
\begin{claim}\label{c:1}
  If $a\in \mathcal A$ is such that $\text{supp}(\pi)\cap a$ is
  infinite then $\text{supp}(\pi)\cap a$ is a co-finite subset of $a$.
\end{claim}
\begin{proof}
 Let $a'\in {\mathcal A}\setminus \{a\}$
If the claim fails then there are
 infinitely many $n \notin a$ such that
 $n \in \text{supp}(\pi)$ but $\pi_{a,a'}(n)\in \text{supp}(\pi)$. For any such
 $n$ it follows that 
$ \pi\circ\pi_{a,a'}(n) \neq \pi_{a,a'}(n)$
while $\pi_{a,a'}\circ\pi(n) = \pi_{a,a'}(n)$.
Hence $[pi]_{\Fin}$ and $[pi_{a,a'}]_{\Fin}$ do not commute.
\end{proof}

\begin{claim}\label{c:2}
  If  $a\in \mathcal A$ is such that $\text{supp}(\pi)\cap a$ is
  infinite then 
$\pi(a) \subseteq^* a$. 
\end{claim}
\begin{proof}
  If not, there are
 infinitely many $n \in a$ such that
  $\pi(n)\notin a$. Let $X$ be
 the set of all such $n$ and choose $a' \in \mathcal A$ such that
 $\pi(X)\setminus a'$ is infinite. 
Then $\pi_{a,a'}\circ\pi(n) = \pi(n)$
and $\pi_{a,a'}(n)\neq n$ for any $n \in \pi^{-1}(X\setminus a')$.
\relax From the last inequality it follows that $\pi(\pi_{a,a'}(n))\neq
 \pi(n)$ and hence 
$[\pi_{a,a'}]_{\Fin}$ does not commute with $[\pi]_{\Fin}$.
\end{proof}
\begin{claim}\label{c:3}
  If  $a\in \mathcal A$ is such that $\text{supp}(\pi)\cap a$ is
  infinite then there is some  $i \in \Integers$ such that
$\pi\restriction a \equiv^* \pi_a^i
  \restriction a$. 
\end{claim}
\begin{proof}
 From Claim~\ref{c:1} and  Claim~\ref{c:2}
it follows that for almost all $n \in a$ 
there is some $k(n)$ such that $\pi(n) =
  \pi_a^{k(n)}(n)$. If the claim is false then there are infinitely many $n\in
  a$ such that  $k(n)\neq k(\pi_a(n))$. For any such $n$ it follows that
$$\pi_a\circ\pi(n) = \pi\circ\pi_a^{k(n)}(n) = \pi_a^{k(n) + 1}(n)\neq
\pi_a^{k(\pi_a(n)) + 1}(n) =
\pi_a^{k(\pi_a(n))}(\pi_a(n)) =
\pi\circ\pi_a(n)$$ and hence
$\pi_a\circ\pi$ and $\pi\circ\pi_a$ disagree on infinitely many integers. 
\end{proof}
There are now two cases to consider.
\begin{case1}
There is a finite subset
$\{a_1, a_2, \ldots a_n\}\subseteq \mathcal A$ such that $\text{supp}(\pi)
\subseteq^* \bigcup_{i=1}^na_i$.   
\end{case1}
In this case, use Claim~\ref{c:3} to choose integers $k_i \in \Integers$
such that $$\pi\restriction a_i
\equiv^* \pi_{a_i}^{k_i}$$ for each $i \leq n$.
This contradicts that $[\pi]_{\Fin} \notin\Phi(F({\mathcal A})) $. 
\begin{case2}
There is no finite subset
$\{a_1, a_2, \ldots a_n\}\subseteq \mathcal A$ such that $\text{supp}(\pi)
\subseteq^* \bigcup_{i=1}^na_i$.   
\end{case2}
In this case there are uncountably many $a\in \mathcal A$ such that 
$\text{supp}(\pi)
\cap a$ is infinite. Use Claim~\ref{c:3} to conclude that
 there is some $i \in \Integers$ such that, without loss of generality,
$\pi\restriction a \equiv^* \pi_a^i$ 
for uncountably many $a\in \mathcal A$.
Hence there is some $k\in \Naturals$ such that
$\pi\restriction \{n\in a : n \geq k\} = \pi_a^i
  \restriction \{n\in a : n \geq k\}$ for uncountably many $a\in \mathcal A$.
Hence there are distinct $a$ and $b$ in $\mathcal A$ such that
it is possible to choose $j \in \{n\in a : n \geq k\}\cap \{n\in b : n
\geq k\}$ such that $\pi(j) = \pi_a^i(j) \neq \pi_b^i(j) = \pi(j)$.
\end{proof}
\section{A lower bound}
The next series or preliminary lemmas will be used in the proof of
Theorem~\ref{t:pt} 
which establishes a lower bound for
$A([\Naturals]^{<\aleph_0})$. Corollary~\ref{c:ps} has as a trivial
consequence the fact  that any
maximal Abelian subgroup of the 
full symmetric group of the integers has cardinality $2^{\aleph_0}$;
however, this can also be shown by using the topology of pointwise
convergence on this group and noting that any maximal Abelian subgroup
must be closed, and hence have cardinality $2^{\aleph_0}$.

\begin{lemma}\label{l:fg}
  Let ${\mathcal S}$ be a finite subset of $\Sacks$ whose elements almost
  commute with each other.
  \begin{enumerate}
  \item If all the orbits of
  each $\pi\in {\mathcal S}$ are finite then
each element of $\Omega_{{\mathcal S}}$ is finite.
 \item If, in addition, for each $\pi \in {\mathcal S}$ all the orbits of
  $\pi$ have size less than or equal to $m(\pi)$ then 
 the cardinality of all but finitely many elements of $\Omega_{{\mathcal S}}$ 
is no greater than $\prod_{\pi \in {\mathcal S}}m(\pi)$.
  \end{enumerate}
\end{lemma}
\begin{proof}
Proceed by induction on $n = |{\mathcal S}|$, the case $n=1$ being trivial. 
If the lemma is true
   for $n$
let ${\mathcal S} = \{\pi_i\}_{i=1}^{n+1}$ and
let ${\mathcal S}' = \{\pi_i\}_{i=1}^n$.
 Define 
$$B = 
   \bigcup_{i=1}^{n}\bigcup_{j=i}^{n+1}\text{NC}(\pi_i,\pi_j)$$ and,
if the orbits of each $\pi \in {\mathcal S}$
are bounded by $m(\pi)$ then 
let $B'$ be the union of those finitely many 
$A\in\Omega_{{\mathcal S}'}$ whose cardinality is not bounded by
  $\prod_{i=1}^{n}m(\pi_i)$.
 Define 
 $$B^* = \Omega_{{\mathcal
      S}'}\left(\text{orb}_{\pi_{n+1}}(B\cup B')\right).$$
Observe that $B^*$ is finite by the induction hypothesis and the fact
the orbits of $\pi_{n+1} $ are finite.
 Hence, it suffices to show
 that if  $C\in \Omega_{{\mathcal S}'}$ and $C \cap B^* =
 \emptyset$
then $C' =  \text{orb}_{\pi_{n+1}}(C)$ belongs to
$\Omega_{{\mathcal S}}$. The fact that it   is finite is immediate
from the hypothesis that all orbits of
  $\pi_{n+1}$ are finite; similarly,
  if  $\card{C} \leq \prod_{i=1}^{n}m(\pi_i)$ then it follows that 
 $\card{C'} \leq \prod_{i=1}^{n+1}m(\pi_i)$. 

To see that $C' \in \Omega_{{\mathcal S}}$
it suffices to show that if $ i \leq n$ and $c \in C'$ then
$\text{orb}_{\pi_i}(c)\subseteq C'$. 
If not then there is some $d \in C$ such
  that $c \in \text{orb}_{\pi_{n+1}}(d)$. Since
  $\text{orb}_{\pi_i}(d)\subseteq C\subseteq C'$ it follows that there must be
  some $e \in  \text{orb}_{\pi_{n+1}}(d)$ such that $\pi_i(e) \in C'$ and
  $\pi_i(\pi_{n+1}(e)) \notin C'$. But $\pi_{n+1}(\pi_i(e)) \in C'$ by
  definition. Hence $\pi_{n+1}\circ\pi_i(e) \neq \pi_i\circ\pi_{n+1}(e)$
  contradicting that $e \notin B$.
\end{proof}

\begin{lemma}\label{l:comm+}
  Let ${\mathcal S} \subseteq \Sacks$ be finite and suppose that $\pi
  \in \Sacks$ and $\theta\in 
  \Sacks$ almost commute with each member of $S$. Then there is a finite set
$Y$ such that
if $\pi\restriction X\cup Y = \theta
  \restriction X\cup Y $ 
then $\pi\restriction \Omega_{\mathcal S}({X}) = \theta
  \restriction \Omega_{\mathcal S}({X})$. Moreover, if
$\pi$ and $\theta$ actually commute with each member of $S$ then $Y$
  can be taken to be the empty set.
\end{lemma}
\begin{proof}
Let $Y' = \bigcup_{\sigma\in {\mathcal S}}
\text{NC}(\sigma,\pi) \cup 
\text{NC}(\sigma,\theta)$ and let $Y = \bigcup_{\sigma\in {\mathcal
    S}}\sigma(Y')\cup \sigma^{-1}(Y')$. 
  Note that  $\Omega_{\mathcal S}(X\cup Y) = \bigcup_{i=0}^\infty
  X^{(i)}$ where $X^{(0)} = X\cup Y$ and 
  $X^{(n+1)} = \bigcup_{\sigma\in {\mathcal S}}
\text{orb}_{\sigma}(X^{(n)})$ and, hence, it suffices to show by induction
  that  
$\pi\restriction X^{(n)} = \theta
  \restriction X^{(n)}$ for each $n$ assuming that
$\pi\restriction X^{(0)} = \theta
  \restriction X^{(0)}$
. To this end, suppose that
$\pi\restriction X^{(n)} = \theta
  \restriction X^{(n)}$ and $x\in X^{(n+1)}$. Then there is some $\bar{x}\in
  X^{(n)}$ and $\sigma \in \mathcal S$ such that $x \in
  \text{orb}_{\sigma}(\bar{x})$. But $\pi(\bar{x}) =\theta(\bar{x})$ and hence
$\sigma^k(\pi(\bar{x})) =\sigma^k(\theta(\bar{x}))$ for any $k$. If $n
> 1$ then $\bar{x}\notin Y$ and it follows that
$\pi(\sigma^k(\bar{x})) =\theta(\sigma^k(\bar{x}))$ for all $k$. Since $x =
  \sigma^k(\bar{x})$ for some $k$ the result follows.

If $n = 1$ it will be shown by induction on $|k|$ that if ${x}\in
  X^{(1)}$ and
$\bar{x}\in
  X^{(0)}$ and $\sigma \in \mathcal S$ are such that $x
  =\sigma^k(\bar{x})$
then $\theta(x) = \pi(x)$.  If $|k| = 0$ this is immediate.
First assume that $k > 0$ and $\theta(\sigma^{k-1}(\bar{x})) = 
\pi(\sigma^{k-1}(\bar{x}))$. If $\sigma^{k-1}(\bar{x}) \notin Y'$ then
$\sigma^{k-1}(\bar{x}) \notin \text{NC}(\pi,\sigma)\cup \text{NC}(\theta,\sigma)$
  and so
$$\theta(\sigma^{k}(\bar{x})) = \sigma(\theta(\sigma^{k-1}(\bar{x})))
=\sigma(\pi(\sigma^{k-1}(\bar{x}))) = \pi(\sigma^{k}(\bar{x}))$$
as required. On the other hand, if $\sigma^{k-1}(\bar{x}) \in Y'$ then
$\sigma(\sigma^{k-1}(\bar{x})) \in Y$ and so
$\theta(\sigma^{k}(\bar{x})) = \pi(\sigma^{k}(\bar{x}))$ in this case
also.
The case that $ k < 0$ is handled similarly. 
\end{proof}

\begin{defin}
  If $H\subseteq \Sacks$ is a subgroup then define $H$ to be {\em strongly
  almost Abelian} if and only  if for each $h\in H$ there is a finite set
  $F(h)\subseteq \Naturals$ such that if and $h_1$ and $h_2$ belong to $H$ 
then $\text{NC}(h_1,h_2)
  \subseteq F(h_1) \cup F(h_2)$.   
\end{defin}
\begin{lemma}\label{l:ps}
  If $H\subseteq \Sacks$ is an uncountable subgroup
 and $F:H \to [\Naturals]^{<\aleph_0}$ attests to the fact that $H$
  is strongly almost Abelian then
there is a perfect set $P\subseteq \Sacks$ and a finite $W\subseteq \Naturals$
such that:
\begin{itemize}
\item There is some $g^* \in H$ such that for all $n\in \Naturals \setminus W$
  and $\pi \in P$ either $\pi(n) = n$ or $\pi(n) = g^*(n)$.
\item   $\text{NC}(\pi,h) \subseteq W \cup F(h)\cup h^{-1}(W)$ for $\pi \in P$ and $h\in H$.
\end{itemize}
Moreover if  $H$ is actually Abelian and not just strongly almost Abelian
 then $W$ can be assumed to be empty and
it can be concluded that each $\pi \in P$ commutes with each $h \in H$.
\end{lemma}
\begin{proof}
   Given $X \subseteq \Naturals$ and a finite  $W \subseteq \Naturals$ define
$cl^0_{W}(X) = X$, 
$$cl^1_{W}(X) = 
\{z \in \Naturals\setminus W: (\exists h \in H)(\exists x \in X\setminus
F(h))
 z = h(x) \AND z\notin F(h^{-1})\}$$
 and let
 $cl^{n+1}_{W}(X) = cl^1_{W}(cl^n_{W}(X))$ and then, let
$cl_{W}(X) = \bigcup_{i=1}^\infty cl^i_{W}(X)$.
Observe first that it follows from an argument similar  to that in
  Lemma~\ref{l:comm+} that,
if $F(g_1) \subseteq W$ and $F(g_2) \subseteq W$
and $g_1\restriction X = g_2\restriction X$ then
 $g_1\restriction cl^1_{W}(X) =g_2\restriction cl^1_{W}(X) $ and hence,
 $g_1\restriction cl_{W}(X) =g_2\restriction cl_{W}(X) $.

If,  for every $W \in [\Naturals]^{<\aleph_0}$ there is some 
$A_W\in [\Naturals]^{<\aleph_0}$ such that 
$A_W \cup cl_{W}(A_W) = \Naturals$ then
it follows that each  $g\in H$ is  determined by its 
 values on $F(g) \cup A_{F(g)}$. This contradicts that $H$ is uncountable.

Therefore it must be the case that
   there is some $W \in [\Naturals]^{<\aleph_0}$ such that
$cl_{W}(A) \neq \Naturals$
  for every $A \in [\Naturals]^{<\aleph_0}$.
Hence, it is possible to choose some $W'\supseteq W$ such that the set
of all $g\in H$ such that $F(g) \subseteq W'$ is uncountable. 
Observe that  $cl_{W}(A)\supseteq cl_{W'}(A)$ for any $A$,
   so it is possible to select
$\{a_i\}_{i=1}^\infty \subseteq \Naturals$ be such that
$\{cl_{W'}(\{a_i\})\}_{i=1}^\infty$ is an infinite family. Observe that
if $g\in H$ is such that $F(g) \subseteq W'$ and $g\restriction
cl_{W'}(\{a_i\})$ is the identity for all but finitely many $i$ then $g$ is
determined by its values on $F(g) \cup \{a_i : (\exists n \in
cl_{W'}(\{a_i\})) g(n) \neq n\}$. Hence, there must be some
  $\bar{g}\in H$  such that $F(\bar{g}) \subseteq W'$ and $\bar{g}\restriction
cl_{W'}(\{a_i\})$ is not the identity for infinitely many $i$.
Let $$Z = \{i \in \Naturals: (\exists n \in
cl_{W'}(\{a_i\})) \bar{g}(n) \neq n \AND \bar{g}^{-1}(W')\cap
   cl_{W'}(\{a_i\}) = \emptyset\}.$$ 
First notice  that it follows from  the definition of $cl_{W'}$ and
   the inclusion
$F(\bar{g})\subseteq W'$ that $\bar{g}(cl_{W'}(\{a_i\}))\subseteq
cl_{W'}(\{a_i\})\cup W'$. 
Hence $\bar{g}\restriction cl_{W'}(\{a_i\})$ is a permutation of
   $cl_{W'}(\{a_i\})$ for each $ i\in Z$.
 Therefore, if for each $t:Z \to 2$ the function $g_t$ is defined
$$g_t(n) =
\begin{cases}
  \bar{g}(n) & \text{if } \ n \in cl_{H,W'}(\{a_i\}) \text{ and } \ t(i) = 0\\
n & \text{otherwise.} 
\end{cases}$$
then $g_t$ is a permutation of $\Naturals$.
It is routine to check that each $g_t(h(n)) = h(g_t(n))$ provided that $n\notin
W'\cup F(h)\cup h^{-1}(W')$.
\end{proof}
\begin{corol}\label{c:ps}
  If $H\subseteq \Sacks$ is an uncountable, maximal
 strongly almost Abelian subgroup then
$\card{H} = 2^{\aleph_0}$. 
\end{corol}
\begin{proof}
The maximality of $H$ implies that it must contain the perfect set of the
conclusion of Lemma~\ref{l:ps}.  
\end{proof}

\begin{lemma}\label{l:iso+}
If $H$ is a maximal Abelian subgroup of $\Sacks/\Fin$ and there are
$$\{[\pi_1]_{\Fin},[\pi_2]_{\Fin},\ldots [\pi_n]_{\Fin}\}\subseteq H$$ 
such
that,
letting ${\mathcal S} = \{\pi_1,\pi_2,\ldots \pi_n\}$,
there are infinitely many different cardinalities of  equivalence classes of
$\equiv_{{\mathcal S}}$, then $\card{H} = 2^{\aleph_0}$.
\end{lemma}
\begin{proof}
  Let $A_j = \cup(\Omega_{\mathcal
  S}\cap [\Naturals]^j)$ and note that
  $\{A_j\}_{j=0}^\infty$ is an infinite set. 
  For each $j \geq 2$ choose some $i \leq n$ such that
$\pi_i\restriction A_j $ is different from the
identity on an infinite subset of $A_j $.  
For each $F:\Naturals \to 2$ define
$$\theta_F(k)= 
\begin{cases}
  \pi_j(k) & \IF k\in A_j\AND F(j) = 1\\
  n & \text{ otherwise}\\
\end{cases}$$
and observe that $\theta_F$ is a bijection. It suffices to show that
$\text{NC}(\theta_F,\pi)$ is finite for each
$\pi \in H$. 

To see that this is so, let $\pi \in H$ and let $j$ be so large that
$$\left(\bigcup_{i=1}^n\text{NC}(\pi_i,\pi)\right)\cap
\left(\bigcup_{i=1}^\infty A_i\right)\subseteq \bigcup_{i=1}^j A_i . 
$$
Hence, if $k \geq j$ then $\pi\restriction A_k$ commutes with
$\pi_i\restriction A_k$ for $i \leq n$. From this it follows that
$\pi\restriction A_k$ is a permutation of $A_k$ because, if $\pi(a) \notin
A_k$  then
 $\pi$ is an ${\mathcal
    S}$-isomorphism from $\Omega_{{\mathcal
    S}}(a)$ onto the $\Omega_{{\mathcal
    S}}(\pi(a))$. This contradicts that 
 $|\Omega_{{\mathcal
    S}}(a)| = k \neq | \Omega_{{\mathcal
    S}}(\pi(a))|$.
\end{proof}
\begin{defin}
  If $g \in \Sacks$ then define $I(g) = \bigcup \{\text{orb}_g(n) :
  \card{\text{orb}_g(n)} = \aleph_0\}$. For a finite set ${\mathcal
  S}\subseteq \Sacks$ define $I^*({\mathcal
  S}) = 
\Omega_{\mathcal
  S}(\bigcup_{\sigma \in {\mathcal
  S}}I(\sigma))$.
\end{defin}
\begin{lemma}\label{l:psu}
  If $H \subseteq \Sacks$ is an uncountable,
 maximal, almost commuting subgroup of size less than
  $2^{\aleph_0}$ then $[\Naturals]^{<\aleph_0}\cup\{I^*({\mathcal
  S})\}_{{\mathcal
  S}\in [H]^{<\aleph_0}}$ generates
  a proper ideal. 
\end{lemma}
\begin{proof}
  If not, let $B\subseteq H$ and $C\subseteq \Naturals$ be finite sets such
  that $I^*(B) \cup C = \Naturals$. Without loss of generality
  it may be assumed that $\text{NC}(b,b')\subseteq C$ for each $b$ and $b'$ in $B$.
  Let $S = \{A \in \Omega_B : A\cap C = \emptyset\}$.
 Observe that each set in $S$ is infinite since it must intersect
  some $I(b)$ where $b\in B$. Moreover, $S$ itself is an infinite set since
  Lemma~\ref{l:comm+} would imply that $H$ is countable otherwise.

  Now let $\Sacks_S$ be 
the symmetric group on $S$ and define $\Phi: H \to
  \Sacks_S$ by $\Phi(h)(s) = t$ if and only if $h(s) \equiv^*
  t$. Observe that $\Phi$ is well defined. To see this suppose
that $s\in S$ and
  $h\in H$ and there are distinct $t$ and $t'$ in $S$
  such that $|h(s)\cap t| = |h(s)\cap t| = \aleph_0$. Then there exist
$i$ and $j$ in $s$ such that
$h(b(i)) = b(h(i))$ for all $b\in B$ and
$h(i) \in t$ and $h(j) \in t'$.
But then, since $\{i,j\} \subseteq s \in \Omega_{B}$,
 there is some $g$ in the subgroup generated by $B$ such that $g(i) =
 j$.
Hence $h(j) = h(g(i)) = g(h(i))$. Furthermore, $h(i) \in t$ and $g$ in
the subgroup generated by $B$ together imply that $g(h(i)) \in t$. However,
$h(j) \in t'$ so  $h(g(i)) \neq g(h(i))$ contradicting the choice of
$i$. A similar argument shows that $\Phi$ is a homomorphism.

Moreover, its image $\Phi(H)$ is an Abelian
  subgroup of $\Sacks_S$.   
  To see this, let $s \in S$. 
If $\Phi(g)(\Phi(h)(s)) \neq \Phi(h)(\Phi(g)(s))$ then,
$g(h(s)) \not\equiv^* h(g(s))$ and hence there are infinitely many $i
  \in s$ such that $g(h(i)) \neq h(g(i))$ contradicting that $h$
  almost commutes with $g$.

  To begin it will be shown that there cannot be 
a perfect set $P\subseteq \Sacks_S$
such that:
\begin{enumerate}
\item \label{lla:1} There is some $g^* \in H$ such that for all $s\in S$
  and $\pi \in P$ either $\pi(s) = s$ or $\pi(s) = \Phi(g^*)(s)$.
\item \label{lla:2}   Every element of $ P$ commutes with every element of $h\in H$. 
\end{enumerate}
To see this suppose that $P$ and $g^*$ contradict the assertion.
For $\pi \in P$ define $\pi^* \in \Sacks$ by
$$\pi^*(i) =
  \begin{cases}
    i & \IF  i\in s\in S \AND \pi(s) = s\\
    g^*(i) & \IF i\in s\in S \AND \pi(s) \neq s
  \end{cases}$$
It suffices to show that $\pi^*$ almost commutes with each $h \in H$.
To see that this is so let $i \in \Naturals\setminus \text{NC}(g^*,h)$
 and let $s\in S$ be such that $i\in s$.
If $\pi(s) = s$ then $h(\pi^*(i)) = h(i) = \pi(h(i))$ the last
equality holding because $h(i) \in \Phi(h)(s)$ and  $\pi(\Phi(h)(s)) = \Phi(h)
(\pi(s)) = \Phi(h)(s)$.
On the other hand,
if $\pi(s) \neq s$ then $h(\pi^*(i)) = h(g^*(i)) = g^*(h(i)) =
\pi^*(h(i))$ 
the last
equality holding because $h(i) \in \Phi(h)(s)$  and
$\pi(\Phi(h)(s)) = \Phi(h)(\pi(s)) \neq \Phi(h)(s)$.

  To see that $\Phi(H)$ is not countable suppose otherwise. 
  To begin, notice  that
there must be some  $A\in\Omega_{\Phi(H)}$
such that $\{h\restriction (\cup A)\}_{h\in H}$ is uncountable ---
keep in mind that $A \subseteq \Omega_B$.
This so because  if  not, then it is easy to find $P$ and
  $g^*$ satisfying conditions~\ref{lla:1}
and~\ref{lla:2}. Simply let $g^*\in H$ be any permutation
$\Phi(g^*)\restriction A$ which is different from the identity
on infinitely many sets in $\Omega_{\Phi(H)}$. Then let $P$ be the set
  if all $g\in \Sacks_S$ such that for all $A\in\Omega_{\Phi(H)}$
  either
$g\restriction A = \Phi(g^*)\restriction A$ or else
$g\restriction A$ is the identity.
If no such $g^*$ exists then it follows that $H$ is countable because
$\{h\restriction (\cup A)\}_{h\in H}$ is countable for each $A
  \in\Omega_{\Phi(H)}$ and each $h$ is the identity on all but
  finitely many   $A
  \in\Omega_{\Phi(H)}$.

Hence, there must be some  $A\in\Omega_{\Phi(H)}$
such that $\{h\restriction (\cup A)\}_{h\in H}$ is uncountable and
hence
there is some $h^*\in H$ such that
 $\{h\restriction \cup A : \Phi(h)\restriction A = \Phi(h^*)\restriction A
\}$ is uncountable. Observe that
if $\Phi(h)\restriction A = \Phi(h')\restriction A$ and there is some
$i \in \cup A $ such that $h(i) = h'(i)$ then
$h\restriction \cup A \equiv^* h'\restriction \cup A$. To see this note
first  that if $\{i,j\}\subseteq s \in A$ then there is 
some $b$ in the group generated by $B$ such that $b(i) = j$.
Since $s \notin C$ it follows that $h(j) = h(b(i)) = b(h(i))=
b(h'(i)) = h'(j)$. 
If $j\in \cup A$ then there are $s_i$ and $s_j$ in $A$ such that
$i \in s_i$ and $j \in s_j$ and there is $\bar{h}\in H$ such that
$\Phi(\bar{h})(s_i) = s_j$.
Since $s_i$ is infinite, there is some $i^* \in s_i\setminus
(\text{NC}(\bar{h},h) \cup\text{NC}(\bar{h},h')) $ such that
$\bar{h}(i^*) \in s_j$.
  Hence $h(\bar{h}(i^*)) = \bar{h}(h(i^*)) = 
\bar{h}(h'(i^*)) = h'(\bar{h}(i^*))$ 
 and, since $\{\bar{h}(i^*),
j\}\subseteq s_j$ it follows that
$h(j) = h'(j)$.
But, since $\{h\restriction (\cup A) : \Phi(h)\restriction A = \Phi(h^*)\restriction A
\}$ is uncountable it is possible to find
$h$ and $h'$ such that there are $i$ and $j$ in $\cup A$ such that
$h(i) = h'(i)$, $h(j) \neq h'(j)$ and $\Phi(h)\restriction A =
\Phi(h^*)\restriction A$.

Since $\Phi(H)$ is uncountable and Abelian 
it  follows from Lemma~\ref{l:ps} that
there exist $P$ and $g^*$ satisfying the conditions~\ref{lla:1}
and~\ref{lla:2}.
  
\end{proof}
The following alternate characterization, due to M. Bell,
 of the cardinal invariant
 $\mathfrak p$ of Definition~\ref{d:d12} will be used in the proof of
 Theorem ~\ref{t:pt}.
 \begin{theor}
The cardinal $\mathfrak p$ is the least cardinal such that there is
a $\sigma$-centered
 partially ordered set\footnote{
A
 partially ordered set is said to be $\sigma$-centred
if it is the union of countably many   
 centred  subsets --- in other words, it is the union
of countably many subsets which
   contain a lower bound for any two of their elements.
} $\Poset$
 and a collection $D$ of  $\mathfrak p$
   dense subsets of $\Poset$ for which there is no centred subset
  $G\subseteq \Poset$ intersecting each member of $D$.  
 \end{theor}
 \begin{proof}
   See \cite{MR83e:03077}.
 \end{proof}

\begin{theor}\label{t:pt}
  If $H \subseteq \Sacks/\Fin$ is an uncountable, maximal Abelian
subgroup then $\card{H} \geq \mathfrak p$ --- in other words, $A([\Naturals]^{<\aleph_0}) \geq \mathfrak p$.  
\end{theor}
\begin{proof}
    If $H \subseteq \Sacks/\Fin$ is an uncountable, maximal Abelian
    subgroup and $\card{H} < \mathfrak p$ then it follows from
    Lemma~\ref{l:psu} that $\{I^*({\mathcal
  S}) : \{[\sigma] : \sigma \in {\mathcal
  S}\} \in [H]^{<\aleph_0}\}$ generates a proper ideal.

    Let $\Poset$ be the partial order consist of all $p = (h^p,{\mathcal S}^p)$
    such that:
    \begin{enumerate}
    \item $h^p$ is a finite involution\footnote{In other words,
        $h^p$ is it own inverse.} 
    \item ${\mathcal S}^p  $ is a finite subset such that if $\sigma\in
      {\mathcal S}^p$ then $[\sigma] \in H$ 
    \end{enumerate}
    and define $p \leq q$ if and only if
    \begin{enumerate}
    \item $h^p \supseteq h^q$
    \item ${\mathcal S}^p \supseteq {\mathcal S}^q$
      \item the domain of $h^p\setminus h^q$ is disjoint from $I^*({\mathcal
  S}^q)$ 
        \item if $j$ is in the domain of $h^p\setminus h^q$ 
and $\sigma \in {\mathcal S}^q$ then 
$\sigma(j)$ is in the domain of $h^p\setminus h^q$ and 
$\sigma(h^p(j)) = h^p(\sigma(j)$. 
    \end{enumerate}
   
It is clear that $\Poset$ is $\sigma$-centred. Moreover, the sets $D_\pi  = \{p
\in \Poset : \pi \in {\mathcal S}^p\}$  are all dense. Furthermore, so are the sets
$$E_n = \{ p\in \Poset : n\in \dom(h^p) \cup I^*({\mathcal
  S}^p)\} .$$ 
To see that this is so, let $p\in \Poset$ be given and suppose that $n\notin
I^*({\mathcal S}^p)$. This implies that the $\equiv_{I^*({\mathcal
  S}^p)}$-equivalence class of $n$ is finite by Lemma~\ref{l:fg}. Now let
$h^q$ be the union of $h^p$ and the identity on the $\equiv_{I^*({\mathcal
  S}^p)}$-equivalence class of $n$ and let $q = (h^q,{\mathcal
  S}^p)$. Then $q \in E_n$ and $q \leq p$.

 Hence, if
$\card{H} < \mathfrak p$ then there is a filter
 $G\subseteq \Poset$ meeting each $D_\pi$ for $\pi \in H$ and $E_n$ for $n\in
 \Naturals$. Define $\pi_G:
 \Naturals \to \Naturals $ by
$$\pi_G(j) = 
\begin{cases}
  h^p(j) & \text{if } \  (\exists p \in G) j \in \dom(h^p) \\ 
  j & \text{if } \  (\forall p \in G) j \notin \dom(h^p)  
\end{cases}
$$
It is easily verified that $\pi_G \in \Sacks$. To see that $\pi_G$ almost
commutes with each member of $H$ let $\pi\in H$. Let $p\in G$ be such that
$\pi \in {\mathcal S}^p$. Then if $j \in \Naturals \setminus \dom(h^p$ there
are two possibilities. If there is some $q \in G$ such that $j$ belongs to the
domain of $h^q$ it is clear that $\pi(\pi_G(j)) = \pi(h^q(j)) = h^q(\pi(j)) =
\pi_G(pi(j))$. However, if there is no $q \in G$ such that $j$ belongs to the
domain of $h^q$ then, by virtue of the fact that $E_j\cap G\neq \emptyset$,
 it must be the case that there is some $q\in G$ such that $j \in I^*(
{\mathcal S}^q)$. Since $\pi \in {\mathcal
  S}^q$ it follows that $\pi(j) \in {\mathcal
  S}^q$. Hence $\pi_G(j) = j$ and $\pi_G(\pi(j)) = \pi(j)$ and so
$\pi(\pi_G(j)) = \pi(h^q(j)) = h^q(\pi(j)) =
\pi_G(pi(j))$.

All that remains to be shown is that the following sets are dense
$$D_{\pi,k} = \{ p \in \Poset : (\exists j\geq k) h^p(j) \neq \pi(k)\}$$
for $\pi \in H$ and $ k\in \Naturals$. To establish this, let $p\in \Poset$ be
given. By Lemma~\ref{l:fg} it follows that each
 $\equiv_{{\mathcal  S}^p}$-equivalence class which is disjoint from $I^*({\mathcal
  S}^p)$ is finite. Moreover, by Lemma~\ref{l:iso+} there must be
infinitely many 
of the same cardinality, and, hence there must be two 
$\equiv_{{\mathcal  S}^p}$-equivalence classes $A$ and $B$
such that:
\begin{enumerate}
\item both $A$ and $B$ are disjoint from $I^*({\mathcal S}^p)$
\item both $A$ and $B$ are disjoint from the domain of $h^p$
\item $k < \min(A)$ 
\item $k < \min(B)$
\item $A$ and $B$ are ${\mathcal
    S}$-isomorphic. 
\end{enumerate}
There are two possibilities. If $\Phi = \pi\restriction A$ then let
  let $h^q$ be the union of $ h^p$ and the  identity on $A$
and let $q = (h^q,{\mathcal
  S}^p)$. 
 Otherwise, let $h^q = h^p \cup \Phi \cup \Phi^{-1}$
and let $q = (h^q,{\mathcal
  S}^p)$. In either case $q \leq p$ and $ q\in D_{\pi,k}$.
\end{proof}
\section{$A([\Naturals]^{<\aleph_0})$ can be smaller than $\mathfrak
  a$}
Through this section the notation $f\leq^* g$ will be used to denote the
relation of eventual domination --- in other words, $f(n) \leq g(n)$
for all but finitely many $n$.
\begin{defin}\label{d:wd}
  A partial order $\Poset$ will be said to be {\em weakly dominated}
  over the model $V$ if 
for every $H:\Naturals \to \Naturals$ belonging to
  $V^{\Poset}$ 
either there is $g: \Naturals \to \Naturals$ in $V$ such that $H\leq^* g$ or
 for every $f:\Naturals \to \Naturals$
  belonging to $V$ 
there is some $R \in \prod_{n=1}^\infty[\Naturals]^{\leq n}$
  belonging to $V$ such that the
  following conditions hold:
  \begin{itemize}
\item $f(\max(\bigcup_{i\in n}R(i)))< \min(R(n))$ 
\item there are infinitely many $n$ such that $H\cap (n\times  R(n))\neq
  \emptyset$.  
  \end{itemize}
\end{defin}

\begin{lemma}\label{l:mm}
  If $G\in \Laver_\alpha$ is generic over $V$, $f\in \fomom$ belongs to $
  V[G]$ and there is $g\in \fomom$ in $V$ such that 
  $f \leq^* g$ then there is $f^*\in\fomom$ in $ V$ such that
$|f\cap f^*| = \aleph_0$.
\end{lemma}
\begin{proof}
  Let $\{n_i\}_{i=0}^\infty$ be an increasing sequence of integers
  such that $n_{i+1} - n_i < n_{i+2} - n_{i + 1}$ for each $i$. In
  $V[G]$ let
the function $\bar{f}$ be defined on $ \Naturals$ by
$\bar{f}(i) = f\restriction [n_i, n_{i+1})$.  As $f$ is bounded by $g
  \in V$ it follows from well known properties of Laver forcing that
there is  $\bar{f}^* \in V$ such that $\bar{f}^*(i) \in
  [\prod_{j=n_i}^{n_{i+1} - 1}g(j)]^{i+1}$ 
and $\bar{f}(i)\in \bar{f}^*(i)$  for each $ i\in \Naturals$.
It is then easy to define $f^*:\Naturals \to \Naturals$ in $V$ such that for
  all $i \in \Naturals$ and $h\in\bar{f}^*(i)$
 there is some $j \in [n_i, n_{i+1})$ such that $f^*(j) = h(j)$.
Hence, $f^*$ is the desired function.
\end{proof}
\begin{lemma}\label{l:star}
  If $\Poset$ is  weakly dominated
  over the model $V$ and $V^{\Poset}$ is a model of \CH then 
there are permutations $\{p_\xi\}_{\xi\in\omega_1}$ of $\Naturals$ in $V$
 which mutually almost commute and which are maximal with respect to
  this property in $V^{\Poset}$.
\end{lemma}
\begin{proof}
  Construct involutions
 $\{p_\xi\}_{\xi\in\omega_1}\subseteq \Sacks$ 
by induction on $\xi$ such
  that any two almost commute.

Using the fact that $V^{\Poset}$ is a model of \CH let
$\{r_\alpha\}_{\alpha\in\omega_1}$ enumerate all $\Poset$-names for
permutations of $\Naturals$ which are forced not to belong to $V$
and suppose that $\{p_\xi\}_{\xi\in\eta}$
have been constructed.

 Now let $\{\eta(i)\}_{i=1}^\infty$ enumerate $\eta$ and define
$$\bar{\Omega}_n =\Omega_{\{p_{\eta(1)}, p_{\eta(2)}, \ldots, p_{\eta(n)}\}} \AND
\bar{\Omega}_n(X) =\Omega_{\{p_{\eta(1)}, p_{\eta(2)}, \ldots,
  p_{\eta(n)}\}}(X)$$
for any set $X$. 
Note that since the $p_\xi$ are almost commuting
involutions it follows that
 $\Omega_n \subseteq^* [\Naturals]^{\leq 2^n}$.
Now, 
for any  $z\in [\Naturals]^{<\aleph_0}$ 
 let $\bar{\tau}_n(z)$ be the structure 
$$\left(|\bar{\Omega}_n(z)|,
  \{\Delta_{\bar{\Omega}_n(z)}\circ
  p_{\eta(i)}\circ \Delta_{\bar{\Omega}_n(z)}^{-1}\}_{i=1}^n\right)$$ 
 Since$\{\tau_n(\{j\})\}_{j\in \Naturals}$ is finite for
each $n$, it follows that, for all $j$, it is possible to choose a finite set
$C(j)\subseteq \Omega_j$ such that 
\begin{equation}
  \label{eq:1}
  \Omega_j( C(j))\cap j = \emptyset
\end{equation}
\begin{equation}
  \label{eq:2}
  \{\tau_j(\{i\}) : i \in C(j)\} = \{\tau_j(\{i\}) :i\in \Naturals \AND \Omega_j(\{i\})\cap j = \emptyset\}
\end{equation}
\begin{equation}
  \label{eq:3x}
(\forall i\in \Naturals)(\forall i'\in C(j))  \text{  if   }
\Omega_j(\{i\})\cap j = \emptyset \AND \tau_j(i)= \tau_j(i')\text{ then }
\min(\Omega_j(\{i'\}) < \min(\Omega_j(\{i\})) . 
\end{equation}

Now define a function  $H$ 
in $V^{\Poset}$ by letting
 $H(j)$ be the least
integer such that 
$H(j)\geq j  $
and
\begin{equation}
  \label{eq:5}
  (\exists i\in C(j))\tau_j(\{i\}) = \tau_j(\{H(j)\})\AND
\Delta_{\bar{\Omega}_j(\{i\})}\circ
  r_{\eta}\circ \Delta_{\bar{\Omega}_j(\{i\})}^{-1}
\neq
\Delta_{\bar{\Omega}_j(\{H(j)\})}\circ
  r_{\eta}\circ \Delta_{\bar{\Omega}_j(\{H(j)\})}^{-1} .
\end{equation}
If no such integer exists then $r_\eta$ is completely determined by
$r_\eta\restriction \bar{\Omega}_j(j\cup C(j))$ 
and hence belongs to $V$ contradicting
the fact that only names forced not to belong to $V$ were enumerated.

In $V$ define $f:\Naturals \to
\Naturals$
as follows. 
First define\footnote{This is not an error ---
 $\max(\bar{\Omega}_k(\{j\}))$ is not intended.}
$L(k) = \max(\bar{\Omega}_k(j))$.
 Then let $f^*(k)$ be large enough that
for each 
$z\in \left[\Naturals\right]^{\leq k}$
there is $\bar{z}\in \left[\Naturals\right]^{\leq k}$ such that
\begin{equation}
  \label{eq:7}
|\bar{\Omega}_k(z)| = |\bar{\Omega}_k(\bar{z})|
\end{equation}
\begin{equation}
  \label{eq:61}
(\forall i\leq k)(\forall m\in z)(\exists \bar{m} \in \bar{z})
\Delta_{\bar{\Omega}_k(z),\bar{\Omega}_k(\bar{z})}\restriction \bar{\Omega}_i(m)\text{ 
  is a bijection onto } \bar{\Omega}_i(\bar{m})
\end{equation}
\begin{equation}
  \label{eq:61a}
(\forall i\leq k)(\forall \bar{m} \in \bar{z})(\exists m\in z)
\Delta_{\bar{\Omega}_k(\bar{z}),\bar{\Omega}_k({z})}\restriction \bar{\Omega}_i(\bar{m})\text{ 
  is a bijection onto } \bar{\Omega}_i({m})
\end{equation}
\begin{equation}
  \label{eq:8}
\bar{z}\subseteq [L(k), f^*(k))  
\end{equation}
 Finally, let $f(k)= L(f^*(k)) = \max(\bar{\Omega}_k(f^*(k)))$.

There are two cases to consider. First assume that $H\leq^* g$ for
some $g \in \fomom\cap V$. In
this case define $a_n$ inductively by setting $a_0 = 0$ and letting
$a_{n+1} = L(f(a_n))$. Using
Lemma~\ref{l:mm}, find $H^* \in V$ 
such that $H(a_n) = H^*(n)$ for 
infinitely many integers $n$. By modifying $H^*$ if necessary, it may be
assumed that $H^*(n) \leq f(a_{n})$ for all $n$. Hence,
$\bar{\Omega}_k(\{H^*(n)\}) \subseteq L(f(a_n)) = a_{n+1}$ for each $n$ and so, by
the definition of $H$, $\bar{\Omega}_k(\{H^*(n)\})
 \subseteq [a_n, a_{n+1})$. Note that if 
$j$ is the unique element of  $C(a_n)$ such that $$\tau_{a_n}(\{j\}) =
\tau_{a_n}\left(\{H^*(n)\}\right)$$ 
then $a_n \leq j < a_{n+1}$ as well.
  Now construct the permutations
$p_\eta^n$ of the interval $[a_n,a_{n+1})$
so that 
$$p_\eta^n \restriction \bar{\Omega}_{a_n}(\{H^*(n)\}) =
\Delta_{\bar{\Omega}_{a_n}(\{H^*(n)\}),\bar{\Omega}_{a_n}(\{j\})}$$
$$p_\eta^n \restriction \bar{\Omega}_{a_n}(\{j\}) =
\Delta_{\bar{\Omega}_{a_n}(\{j\}),\bar{\Omega}_{a_n}(\{H^*(n)\})}$$
and $p_\eta^n$ is the identity elsewhere on $[a_n,a_{n+1}))$.
Finally, let $p_\eta = \bigcup_{i=0}^\infty p_\eta^i$ and note that
this is a well defined permutation.
It is immediate that the composition $p^i_\eta \circ p^i_\eta$ is the
identity and that for each $j \leq i$, $p^i_\eta\circ p_{\eta(j)} =
p_{\eta(j)}\circ p_\eta^i$. Hence $p_\eta$ almost commutes with each
$p_{\eta(j)}$. For any $n$ such that $H^*(n) = H(a_n)$ it
follows that $p_\eta^n$ does not commute with $r_\eta\restriction
\bar{\Omega}_{a_n}(\{H^*(n)\})$. 

In the second case,  suppose that $R\in  \prod_{n=1}^\infty[
  \Naturals]^{\leq n}$ is in $V$ and
witnesses the conditions of Definition~\ref{d:wd} with respect to $H$
and $f$. Let $\bar{R}(i) = \max(\bigcup_{j\in i}R(i))$.
The permutation $p_\eta$ will be defined so that
$p_\eta = \bigcup_{i=0}^\infty p_\eta^i$ where the $p_\eta^i$ are
constructed by induction so that
\begin{enumerate}
\item The domain and range of $p_\eta^i$ are $[L(\bar{R}(i)),
  L(\bar{R}(i+1)))$. 
\item The composition $p^i_\eta \circ p^i_\eta$ is the identity.
\item For each $j \leq i$, $p^i_\eta\circ p_{\eta(j)} =
  p_{\eta(j)}\circ p_\eta^i$. 
\end{enumerate}
To see that this can be done suppose that $\{p_\eta^j\}_{j\in i}$ have
been constructed. Then the domain and range  of $\bigcup_{j\in
  i}p_\eta^j$ are $L(\bar{R}(i))$ by the first induction hypothesis. Hence,
$p^i_\eta$ will be defined so that its domain and range are
$[L(\bar{R}(i)),L(\bar{R}(i+1))$. To see that this can be
done recall that
$R(i) \in [\Naturals]^{\leq i}$
and, furthermore, if $b \in R(i)$ then $b \geq f(\bar{R}(i))$ and so
$\bar{\Omega}_i(\{b\})\cap f^*(\bar{R}(i)) =
\emptyset$. 
The definition of $f^*$
guarantees that there is some 
$\bar{z}\in [[L(\bar{R}(i)),
f^*(\bar{R}(i)))]^{\leq i}$   
such that
$|\bar{\Omega}_i(\bar{z})| = |\bar{\Omega}_i(R(i))|$
and conditions~\ref{eq:61} and \ref{eq:61a} hold.
Observe that $\bar{\Omega}_i(\bar{z})\cap \bar{\Omega}_i(R(i)) =
\emptyset$ 
because $\bar{\Omega}_i(R(i))\cap
f^*(\bar{R}(i)) = \emptyset$ and  
$\bar{\Omega}_i(\bar{z})\subseteq
f^*(\bar{R}(i))$.
It follows that if
  $p_\eta^i $ is defined so that
$$p_\eta^n \restriction \bar{\Omega}_i(R(i)) =
\Delta_{\bar{\Omega}_i(R(i)),\bar{\Omega}_i(\bar{z})}$$
$$p_\eta^n \restriction \bar{\Omega}_i(\bar{z}) =
\Delta_{\bar{\Omega}_i(\bar{z}),\bar{\Omega}_i(R(i))}$$
and $p_\eta^n$ is the identity elsewhere on 
 $[L(\bar{R}(i)), L(\bar{R}(i+1)))$ then $p_\eta^i $ is a
bijection of $[L(\bar{R}(i)), L(\bar{R}(i+1)))$
and $p_\eta^i \circ p_\eta^i $ is the identity.
Let $p_\eta = \bigcup_{i=0}^\infty p_\eta^i$.
 Note that 
$\max({\bar{\Omega}_i(\bar{z})})
< \min(\bar{\Omega}_i(R(i)))$, $i \leq \max(\bar{\Omega}_i(R(i)))$ and
hence $\max(\bar{\Omega}_i(\bar{z}) \cup \bar{\Omega}_i(R(i))) \leq 
L(\bar{R}(i+1))$. 
Therefore the first two
induction hypotheses hold. 
By observing that 
$p_\eta$ is an isomorphism between 
$\tau_i(\bar{\Omega}_i(R(i)))$ and $\tau_i(\bar{\Omega}_i(\bar{z}))$  
it follows that the last induction hypothesis
  holds as well. 

Now suppose  that  $H(a) = b$ and $b \in R(i)$ and $a \in i$.
Then there is some $\bar{b} \in \bar{z}$ such that
$\Delta_{\bar{\Omega}_i(R(i)),\bar{\Omega}_i(\bar{z})}(\Omega_a(\{b\})) = \Omega_a(\{\bar{b}\})$.
Since $\max({\bar{\Omega}_i(\bar{z})})
< \min(\bar{\Omega}_i(R(i)))$, 
it follows that $\max(\Omega_a(\{\bar{b}\}))
< \min(\Omega_a(\{b\}))$. 
 By the fact that $b$ is
the least integer satisfying  \ref{eq:5}
it follows that there is some $j\in C(a)$ such that
 $\tau_a(\{b\}) = \tau_a(\{j\})$ and
$\Delta_{\bar{\Omega}_a(\{j\})}\circ
  r_{\eta}\circ \Delta_{\bar{\Omega}_a(\{j\})}^{-1}
\neq
\Delta_{\bar{\Omega}_a(\{b\})}\circ
  r_{\eta}\circ \Delta_{\bar{\Omega}_a(\{b\})}^{-1}$
and, hence there is some $n$ such that
$$
  r_{\eta}\circ \Delta_{\bar{\Omega}_a(\{j\})}^{-1}(n)
\neq \Delta_{\bar{\Omega}_a(\{j\})}^{-1}(
\Delta_{\bar{\Omega}_a(\{b\})}\circ
  r_{\eta}\circ \Delta_{\bar{\Omega}_a(\{b\})}^{-1}(n)) =
\Delta_{\bar{\Omega}_a(\{b\}),\bar{\Omega}_a(\{j\})}
\circ
  r_{\eta}\circ \Delta_{\bar{\Omega}_a(\{b\})}^{-1}(n))
=$$$$
\Delta_{\bar{\Omega}_a(\{b\}),\bar{\Omega}_a(\{j\})}
\circ
  r_{\eta}\circ
  \Delta_{\bar{\Omega}_a(\{b\})}^{-1}(\Delta_{\bar{\Omega}_a(\{j\})}\circ\Delta_{\bar{\Omega}_a(\{j\})}^{-1}(n)) 
$$
and therefore
 $r_\eta(\Delta_{\bar{\Omega}_a(\{j\})}^{-1}(n))\neq (p^i_\eta\circ r_\eta\circ p^i_\eta) (
 \Delta_{\bar{\Omega}_a(\{j\})}^{-1}(n))$.  
It follows that  $p_\eta\circ
r_\eta $ and $r_\eta\circ
p_\eta$ disagree on infinitely many integers.
\end{proof}

\begin{notat}
  The notation $\Laver$ will be used to denote the Laver partial
 order. Most of the notation and terminology regarding this partial
 order will be taken from \cite{MR54:10019}.
  Let $E_T:\wfomom \to T$ be the unique
  bijection preserving the lexicographic ordering from $\wfomom$ onto
 the nodes above the root of $T$.
If $T\in \Laver$ and $t \in \wfomom$ let $T\langle t \rangle$ denote
$E_T(t)$. If $s\in T$ then $T_s$ will denote
the subtree of $T$ consisting all nodes comparable to $t$.
 Fix an enumeration
 $\{s_n\}_{n=0}^\infty$ of $\wfomom$ such that if $s_n \subseteq s_m$
 then $n \leq m$.
Given $T \in \Laver$ and $ n\in \Naturals$ let $\{S_{n,j}(T)\}_{n=0}^{j}$
 list the components\footnote{See the definition on page~156
  of \cite{MR54:10019}.}  of $T$ determined by $\{s_n\}_{n\in j}$ ---
 to be precise, $S_{n,j}(T)$ is the subtree of $T\langle s_n\rangle$
 consisting of all nodes $t$ such that if $n < m \leq j$ then
 $s_m\not\subseteq t$. 

The notation $\Laver_\alpha$ will be used to denote the countable
  support iteration of $\alpha$ Laver partial orders. For the
  definition of $p\geq^n_F q$ see \cite{MR54:10019}.
\end{notat}

\begin{lemma}\label{l:kl}
  If $p \in \Laver_\alpha$,  $n\in\Naturals$, $\forcingname{m}$ is a
  $\Laver_\alpha$-name for an integer  and 
$F$ is a finite subset of the domain of $p$ 
 then   there is some ${p}'$ such that
 \begin{itemize}
 \item $p'(0)\geq p(0)$
\item $p'(0)\langle 0 \rangle = p(0)\langle 0 \rangle$
\item $p'\restriction [1,\alpha) \geq_F^n p\restriction [1,\alpha)$
\item there is $U\subseteq \Naturals$ such that $|U|\leq (n+1)^{|F|}$
and for every integer $k$, for all but finitely many immediate
successors $t$ of $p(0)\langle 0 \rangle$
$${p}'(0)_t {}^\frown p'\restriction [1,\alpha)
 \forces{\Laver_\alpha}{\forcingname{m}\in U\cup (\Naturals \setminus \check{k})}$$
 \end{itemize}
\end{lemma}
\begin{proof}
 If $n = 0$ then this follows from Lemma~12 of \cite{MR54:10019}.
 Let $\forcingname{a}$ of that Lemma be $ \frac{1}{\forcingname{m}}$ and
 find $u$ and $p'$ such that for all $\epsilon > 0$
$${p}'(0)_t {}^\frown p'\restriction [1,\alpha)
 \forces{\Laver_\alpha}{|\forcingname{a} - u| < \epsilon}$$
for all but finitely
 many immediate successors $t$ of $p(0)\langle 0 \rangle$.
If $u \neq 0$ let $\epsilon$ be small enough that there is some 
$U = \{j\}$ such that $|u - 1/i| >
 \epsilon$ for all $i \geq 1$ unless $i = j$.
Then for any immediate
successor $t$ of $p(0)\langle 0 \rangle$
such that
${p}'(0)_t {}^\frown p'\restriction [1,\alpha)
 \forces{\Laver_\alpha}{|\forcingname{a} - u| < \epsilon}$
it follows that
${p}'(0)_t {}^\frown p'\restriction [1,\alpha)
 \forces{\Laver_\alpha}{\forcingname{m}\in U}$ and hence
$${p}'(0)_t {}^\frown p'\restriction [1,\alpha)
 \forces{\Laver_\alpha}{\forcingname{m}\in U}.$$
 
If $u = 0$ and $k\in \Naturals$ let $\epsilon < 1/k(k+1)$. 
 Then for any immediate
successor $t$ of $p(0)\langle 0 \rangle$
such that
${p}'(0)_t {}^\frown p'\restriction [1,\alpha)
 \forces{\Laver_\alpha}{|\forcingname{a} - u| < \epsilon}$
it follows that
${p}'(0)_t {}^\frown p'\restriction [1,\alpha)
 \forces{\Laver_\alpha}{\forcingname{m}\in \Naturals \setminus \check{k}}$.

The case for $n>0$ follows from the case $n=0$ by induction and a
counting argument.
\end{proof}
\begin{lemma}\label{t:ml}
  For any ordinal $\alpha$ the partial order $\Laver_\alpha$ is weakly
  dominated over $V$. 
\end{lemma}
\begin{proof}
Let  $f : \Naturals \to \Naturals$
  belong to $V$ and  $H$ be an $\Laver_\alpha$-name for a function
from $\Naturals$ to $\Naturals$. Let $p \in \Laver_\alpha$.
The function  $R\in \prod_{i=1}^\infty[ \Naturals]^{\leq i}$
 satisfying the three
  conditions of Definition~\ref{d:wd}
will be constructed by induction in $V$ unless $H$ is dominated by a
  function in $V$.

Let $\theta$, $\beta$ and
$\delta$ be monotone functions from $\Naturals$ to $\Naturals$ such that
\begin{itemize}
\item  $\lim_{n\to\infty}\theta(n)=\lim_{n\to\infty}\delta(n) =
\lim_{n\to\infty}\beta(n) =\infty$
\item if $s_m$ is the predecessor or $s_n$ then $\beta(n)\leq \beta(m)
  + 1$ 
\item $\beta(n)\cdot(\theta(n)+1)^{\delta(n)}\leq n$ 
\item $\beta(n) < n$.
\end{itemize}
Then construct   $\{(p_n,F(n), R^n) \}_{n=0}^\infty$ 
by induction to satisfy the following conditions.
\begin{enumerate}
\item \label{C:B100} $p_0 = p$
\item $p_{n} \leq_{F(n)}^{\theta(n)} p_{n+1}$
\item \label{C:B8} $|F(n)| = \delta(n)$ and $F(n)\subseteq [1,\alpha)$
\item \label{C:B9} $F(n+1)\supseteq F(n) $
\item $p_n(0)\langle s_i\rangle =p_{n+1}(0)\langle s_i\rangle$ if $i
  \leq n$
\item \label{C:B1}if $n\geq 1$ then for each $a\leq \beta(n)$ there is $U^n_a\subseteq
  \Naturals$   such that $$R^{n} = \bigcup_{a\in \beta(n)} U_a^n
\setminus f(\max(R^{n-1}))$$
\item\label{C:B2} $|U^n_a| \leq (\theta(n) + 1)^{\delta(n)}$ for each
  $a\leq \beta(n)$ 
\item \label{C:B4-}
for  all  $m\leq n$ and $a \leq \beta(m)$ and $k \in \Naturals$
 and for all but finitely many
$t$ which are immediate successors of $S_{m,n}(p_n(0))\langle 0 \rangle$
$$S_{m,n}(p_n(0))_t {}^\frown p_n\restriction [1,\alpha)
 \forces{\Laver_\alpha}{H(\check{a})\in
   \check{U}^m_a\cup (\Naturals\setminus 
k
)}$$ 
\item \label{C:B4}
for  all  $m\leq n$ and $a \leq \beta(m)$
$$S_{m,n}(p_n(0)) {}^\frown p_n\restriction [1,\alpha)
 \forces{\Laver_\alpha}{H(\check{a})\in
   \check{U}^m_a\cup (\Naturals\setminus 
f(\max(R^{n}))
)} .$$ \end{enumerate}

It will first be shown that this suffices. Let $p_\omega$ be the limit
of the $p_n$ and define $R:\Naturals\to [ \Naturals]^{\leq \aleph_0}$ by 
letting $R(n)=R^n$ for each integer $n$. 
Let $f^*$ be defined by $$f^*(n) = f\left(\max\left(\bigcup_{\beta(i)
  \leq n}R^i\right)\right).$$  
If $H\leq^* f^*$ then $H$ is dominated by a function in $V$ and
there is nothing to do; so suppose that
$H\not\leq^* f^*$.
\relax From Conditions~\ref{C:B1} and~\ref{C:B2} and the choice
of $\beta$,  $\theta$ and $\delta$ it follows that 
$R\in \prod_{n=1}^\infty [\Naturals]^{\leq n}$.
Furthermore, from Condition~\ref{C:B1} it follows  that
$f(\max(\bigcup_{i\in n}R(i)))< \min(R(n))$ for each $n$. To see that
$p_\omega$ forces the last requirement to be satisfied, suppose not and that
$q\geq p_\omega$ and $K$ are such that 
$q\forces{\Laver_\alpha}{(\forall n \geq K)H\cap (n\times R(n)) = \emptyset}$.
Extend $q$ so that 
$q\forces{\Laver_\alpha}{H(\check{a}) = \check{b} \geq
  \check{f}^*(\check{a})}$ for some integers $a\geq  
K$ and $b$ and
so that $q(0)\langle0\rangle = s_m$ where $\beta(m) > a$ and $b
<f(\max(R^m))$. Let $m'$ be the least integer such that $s_{m'}\prec
s_m$ and $b
<f(\max(R^{m'}))$ and let $m''$ be such that $s_{m''}$ is 
the predecessor of $s_{m'}$. 
It follows that $b
\geq f(\max(R^{m''}))$. Also observe that, if $\beta(m') \leq a$ then
$b \geq f^*(a)$ implies that $b\geq f(\max(R^{m'}))$ which is a
contradiction. Hence, $\beta(m') > a$ and, by the hypothesis on
$\beta$, $\beta(m'') \geq a$. If it could be established that
$b\in U^{m''}_a$ it would follow that $(a,b) \in (\beta(m'') +
1)\times R^{m''}\subseteq m''\times R(m'')$ and this would suffice. 

There are now two cases to consider: Either 
$b <f(\max(R^{m'-1}))$ or $b \geq f(\max(R^{m'-1}))$.
In the first case
use Condition~\ref{C:B4}
to conclude
that
$$S_{m'',m'-1}(p_{m'-1}(0)){}^\frown p_{m'-1}\restriction [1,\alpha)
\forces{\Laver_\alpha}{H(\check{a})\in 
   \check{U}^{m''}_a}.$$  Because
$S_{m'',m'-1}(p_{m'-1}(0)){}^\frown p_{m'-1}\restriction [1,\alpha)
\leq q$ it follows that    
$b$ is not ruled out as a possible value for $H(\check{a})$
by $S_{m'',m'-1}(p_{m'-1}(0)){}^\frown p_{m'-1}\restriction [1,\alpha)$
 and, hence,
 that $b\in U^{m''}_a$ as required.
In the second case, it would suffice to show that
$b\in U^{m'}_a$ because then it would follow that
 $(a,b) \in (\beta(m') +
1)\times R(m')$. But this is clear since
$$S_{m',m'}(p_{m'}(0)){}^\frown p_{m'}\restriction [1,\alpha)\leq q$$ and
so $b$ is not ruled out as a possible value for $H(\check{a})$ in this
case either.

To show that the induction can be carried out, let $p_0 = p$ and
suppose that 
 $\{(p_n,F(n), R^n) \}_{n\in j}$ 
have been constructed satisfying the required induction
hypotheses. Let $F(j)$ be given by some bookkeeping scheme so that
Conditions~\ref{C:B8} and~\ref{C:B9}  are satisfied.
 Using Lemma~\ref{l:kl} find
${p}$ such that
 \begin{itemize}
 \item $p(0)\geq S_{j,j}(p_{j-1}(0))$
\item $p(0)\langle 0 \rangle = S_{j,j}(p_{j-1}(0))\langle 0 \rangle = 
p_{j-1}(0)\langle s_j \rangle$
\item $p\restriction [1,\alpha) \geq_{F(j)}^{\theta(j)}
 p_{j-1}\restriction [1,\alpha)$
\item there are $U^j_a\subseteq \Naturals$ for each $a\leq \beta(j)$
  such that $|U^j_a|\leq (\theta(j)+1)^{|\delta(j)|}$
and for every integer $k$, for all but finitely many immediate
successors $t$ of $p_{j-1}(0)\langle s_j\rangle$
$${p}(0)_t {}^\frown p\restriction [1,\alpha)
 \forces{\Laver_\alpha}{H(a) \in \check{U}^j_a
   \cup(\Naturals
\setminus \check{k})}.$$
 \end{itemize}
Then let $R^j$ be defined according to Condition~\ref{C:B1}.
Now, by removing finitely many immediate successors of
$p_{j-1}(0)\langle s_j\rangle $ from $p(0)$ it is possible to
 obtain $\bar{p}_j\subseteq p(0)$ such that $\bar{p}_j\langle 0
 \rangle =  p(0)\langle 0
 \rangle$ and
 for each $a\leq \beta(j)$
 $$\bar{p}_j {}^\frown p\restriction [1,\alpha)
 \forces{\Laver_\alpha}{H(a) \in \check{U}^j_a
   \cup(\Naturals
\setminus \max(R^j))}.$$
Similarly, but using the induction hypothesis for each $n\leq j-1$
to remove finitely many immediate successors of
$p_{j-1}(0)\langle s_n\rangle $ from $S_{n,j-1}(p_{j-1}(0))$,
 it is possible to
 obtain $\bar{p}_n\subseteq S_{n,j-1}(p_{j-1}(0))$ 
such that $\bar{p}_n\langle 0
 \rangle =  p_{j-1}(0)\langle s_n
 \rangle$ and
 for each $a\leq \beta(n)$
 $$\bar{p}_n {}^\frown p_{j-1}\restriction [1,\alpha)
 \forces{\Laver_\alpha}{H(a) \in \check{U}^n_a
   \cup(\Naturals
\setminus \max(R^j))}.$$
Let $$p_j^* = \bar{p}_j {}^\frown p\restriction [1,\alpha)$$
and, for $n \leq j-1$, let
$$p^*_n = \bar{p}_n {}^\frown p_{j-1}\restriction [1,\alpha)$$
and then define $p_j$ to be the join of $\{p^*_i\}_{i=0}^j$.
It is immediate to check that all of the induction hypotheses
are satisfied by $p_j$, $F(j)$ and $R^j$.
\end{proof}
The following result, due to S. Shelah, is 5.31 in \cite{math.LO/9712283}. It will be
useful to know that Laver forcing $\Laver$ is NEP.
\begin{lemma}\label{5.31-630}
Let $\{B_\alpha\}_{\alpha\in\omega_1}$ be family of Borel sets in the
  model of set theory $V$ such that
  $V\models\bigcap_{\alpha\in\omega_1}B_\alpha=\emptyset$.
Let $\Poset$ be a NEP partial order with definition in $V$ 
  and suppose that $\{\Poset_\alpha\}_{\alpha\in\omega_2}$ is a countable
  support iteration such that
  $\Poset_{\alpha + 1} = \Poset_\alpha * \Poset$ for any $\alpha\in \omega_2$.
If 
$$1\forces{\Poset_{\alpha +
    1}}{\bigcap_{\alpha\in\omega_1}B_\alpha=\emptyset}$$ 
 for each $\alpha\in \omega_1$ then
$$1\forces{\Poset_{\omega_2}}{\bigcap_{\alpha\in\omega_1}B_\alpha=\emptyset}
.$$ 
\end{lemma}
\begin{theor}\label{tt:ml}
  It is consistent that $A([\Naturals]^{<\aleph_0}) = \aleph_1 <
  \mathfrak a$.
\end{theor}
\begin{proof}
The model witnessing this is the one obtained by forcing with
$\Laver_{\omega_2}$ over a model $V$ satisfying \CH . From
Lemma~\ref{l:star} and 
Lemma~\ref{t:ml}  it follows that
there are permutations $\{p_\xi\}_{\xi\in\omega_1}$ of $\Naturals$ in $V$
 which mutually almost commute and which are maximal with respect to
  this property; in other words, letting $B_\alpha$ be the Borel  set of all
  permutations of the integers which almost commute with $p_\alpha$
  but are not equal to $p_\alpha$ modulo a finite set,
$\bigcap_{\alpha \in \omega_1} B_\alpha = \emptyset$.
Moreover, it follows from Lemma~\ref{l:star} and 
Lemma~\ref{t:ml} that
 $$1\forces{\Laver_{\omega_1}}{\bigcap_{\alpha \in \omega_1} B_\alpha
   = \emptyset}$$ and, hence by Lemma~\ref{5.31-630} it follows that
 $$1\forces{\Laver_{\omega_2}}{\bigcap_{\alpha \in \omega_1} B_\alpha
   = \emptyset} $$
or, in other words,
$1\forces{\Laver_{\omega_2}}{A([\Naturals]^{<\aleph_0}) = \aleph_1}$. 

The fact that ${\mathfrak a} = \aleph_2$ in this model is well known
and can be found, for example, in \cite{ba.ju.book}.
\end{proof}

\section{The Cohen model and the summable ideals}
The ideal ${\mathcal I}_{1/x}$ is defined to be the set of all
$X\subseteq \Naturals$ 
such that $\sum_{x\in X}1/x < \infty$.  
It will be shown  that $\Sacks({\mathcal I}_{1/x})/\Fin({\mathcal
    I}_{1/x})$ has a maximal Abelian subgroup of size $\aleph_1$ in any
  model obtained by adding uncountably many Cohen reals for any
  function $h$ such that  $h(n) \geq 1/n$ for all but finitely many
  $n\in \Naturals$.
The basic scheme of the argument is that in a model of the form
$V[\{c_\xi\}_{\xi\in\omega_1}] $ where $\{c_\xi\}_{\xi\in\omega_1}$
are Cohen reals, it is possible to define permutations
$\{\pi_\xi\}_{\xi\in\omega_1}$ such that $\pi_\xi \in
V[\{c_\eta\}_{\eta\in\xi + 1}]$ which form an Abelian subgroup which
is close to maximal in the following sense: If
$G$ is the group generated by $\{\pi_\xi\}_{\xi\in\omega_1}$,
 $G'$ is the group generated by $G$ and all permutations which are
 first order definable from $G$ and $G''$ is a maximal Abelian
 subgroup of $G'$ containing $G$ then $G''$ is actually a maximal
 Abelian subgroup.
Consequently, for most of the rest of  this section a family of
permutations with certain properties will be fixed --- these should
thought of as the permutations  obtained from the Cohen reals. Some
notation will first be established. 
 
Suppose that $\{\pi_\xi\}_{\xi\in\alpha}$ have been constructed. 
Let $\{\alpha_i\}_{i=0}^\infty$ enumerate $\alpha$ and, by re-indexing, 
suppose that  $\pi_n = \pi_{\alpha(n)}$. Let $\Omega_m =
\Omega_{\{\pi_i\}_{i\in n}}$. 
The elements of $\Omega_m$ will be enumerated 
as $\{\Omega_m^i\}_{i=0}^\infty$
in such a way that $i < j$ if and only if $\min(\Omega_m^{i}) <
  \min(\Omega_m^{j}) $.
The following technical definition will be the key to the induction
hypothesis of the construction.
\begin{defin}\label{d:1}
  A family of permutations ${\mathcal F}
  =\{\pi_n\}_{n=0}^\infty$ will be said to be {\em nice}
if
\begin{enumerate}
\item $\lim_{n\to\infty} \frac{\pi_k(n)}{n} = 1$ for each $k$
\item $\pi_k = \pi_k^{-1}$ for each $k$
\item $\pi_k$ and $\pi_m$ almost commute
for each $k$ and $m$ 
\item for every $m\in \Naturals$
there is an integer $K_m$
such that if $i$ and $j$ are greater than $K_m$ then 
$\Omega_m^{i}$ and $\Omega_m^{j}$
are $\{\pi_i\}_{i\in m}$-isomorphic. (See Definition~\ref{d:324})
\end{enumerate}
\end{defin}
\begin{defin}
\label{d:po}
Given a  nice family $\mathcal F$ 
  define the partial order $\Poset({\mathcal F})$ to consist of
  triples $p = (f_p,m_p,\epsilon_p)$ such that
$m_p \in \Naturals$, $\epsilon_p > 0$ and
   there is an integer $I_p$ such that
$f_p$ is a permutation of $\Omega_{m_p}(I_p)$ such that $f_p^{-1} = f_p$ 
and
\begin{equation}
  \label{eq:A0}
 \Omega_{m_p}^{i}\text{ is $\{\pi_i\}_{i\in m}$-isomorphic to }\Omega_{m_p}^{j} \text{ if }
i > I_p \AND j > I_p
\end{equation}
\begin{equation}
  \label{eq:A0+}
(\forall n\in\Naturals \setminus \dom(f_p))
(\forall i \leq m_p)
(\forall j \leq m_p)
 \pi_i\circ\pi_j(n) = \pi_j\circ\pi_i(n)
\end{equation}

and if
 $$\delta = \sup_{j \in \Naturals\setminus \dom(f_p)}\left(\sup_{n \leq m_p}
\left | 1 - \frac{\pi_n(j)}{j} \right |\right)$$
then 
\begin{equation}
  \label{eq:A2.1}
\left(1 +
  2^m\left(
1 - \frac{1}{(1 + \delta)^m}\right)\right)(1 + \delta)^m < 1 + \epsilon_p .
\end{equation}
Define $p\leq q$ 
if  $\epsilon_p \leq \epsilon_q$, $m_p \geq m_q$ and 
\begin{equation}
  \label{eq:A3}
  f_p = f_q \cup \bigcup_{\{u,v\}\in t}\Delta_{\Omega_{m_p}^u,\Omega_{m_p}^v}
\end{equation}
where $t$ partitions $I_p\setminus I_q$ into pairs, and
\begin{equation}
  \label{eq:A4}
  1- \epsilon_q < \frac{f_p(i)}{i} < 1 + \epsilon_q
\end{equation}
for each $i$ in the domain
  of $f_p\setminus f_q$.  
\end{defin}
\begin{lemma}\label{l:cm1}
Let  $\mathcal F$ be a nice family, $m\in\Naturals$ and $K_m$ be as in
Definition~\ref{d:1} and suppose that $i$ and $j$ are both greater
than $K_m$. Suppose further that $k\geq i$ and $x \in \Omega_m^i$ then 
\begin{equation}
  \label{eq:A0x}
(1 - \epsilon) < \frac{\pi_n(x)}{x} < (1 + \epsilon)  
\end{equation}
for each $n \in m$. Then the following inequalities hold:
\begin{equation}
  \label{eq:2c}
  \frac{\max(\Omega_m^{i})}{\min(\Omega_m^{i})} < (1 + \epsilon
  )^m
\end{equation}
\begin{equation}
  \label{eq:3c}
  \frac{\min(\Omega_m^{i+1})}{\min(\Omega_m^{i})} < 1 +
  2^m\left(
1 - \frac{1}{(1 + \epsilon)^m}\right)
\end{equation}
\begin{equation}
  \label{eq:4c}
  \frac{\min(\Omega_m^{i+k})}{\min(\Omega_m^{i})} < \left(1 +
  2^m\left(
1 - \frac{1}{(1 + \epsilon)^m}\right)\right)^k
\end{equation}
and, furthermore, if \ $\Phi : \Omega_m^i \to \Omega_m^{i+k}$ is
a $\{\pi_i\}_{i\in m}$-isomorphism then
\begin{equation}
  \label{eq:5c}
 (1 -  \epsilon)^m < \frac{\Phi(n)}{n} < \left(1 +
  2^m\left(
1 - \frac{1}{(1 + \epsilon)^m}\right)\right)^k(1 + \epsilon)^m
\end{equation}
\end{lemma}
\begin{proof}
In order to prove \ref{eq:2c}  the first thing to note is that if $a \in \Omega_m$ and $\{x,
  y\}\subseteq a$ then there is $k\leq m$ and a sequence $(u_1,u_2,
  \ldots , u_k) \in m^k$ such that 
$$x = \pi_{u_1}\circ\pi_{u_2}\circ\ldots \circ\pi_{u_k}(y) .$$
Given $x \in \Omega_m^{i}$ let $k(x)$ be the least integer such that
$$x = \pi_{u_1}\circ\pi_{u_2}\circ\ldots\circ
\pi_{u_{k(x)}}(\min(\Omega_m^{i}))$$
and proceed by induction on $k(x)$. 
If $k(x) = 0$ then $x = \min(\Omega_m^{i})$ and the result is clear.
Suppose that the lemma has been established for all $x$ such that
$k(x) = n$. Given $x$ such that $k(x) = n+ 1$ it is possible to find
$x'$ such that $k(x') = n$ and $x = \pi_u(x')$ for some $u \in m$.
\relax From \ref{eq:A0x} it follows that $$(1 - \epsilon) 
<\frac{\pi_u(x')}
{x'} <   
(1 + \epsilon)$$ and from the induction hypothesis it follows that 
$$(1 - \epsilon)^n
< \frac{x'}{\min(\Omega_m^{j})} < (1 + \epsilon)^n .$$
Hence,  
$$(1 - \epsilon)^{n + 1}
< \frac{x}{\min(\Omega_m^{j})} < (1 + \epsilon)^{n+1} $$
as desired.

To see that \ref{eq:3c} holds begin by observing that if $i' < i$ and
$\Omega_m^{i'}\setminus \min(\Omega_m^{i}) \neq \emptyset$ then, by
\ref{eq:2c},  $$\min(\Omega_m^{i}) < \max(\Omega_m^{i'}) <
\min(\Omega_m^{i'})(1 + \epsilon)^m$$ 
and hence $$\min(\Omega_m^{i'}) > \frac{\min(\Omega_m^{i})}{(1 + \epsilon)^m}.$$
Therefore, the cardinality of
$$\bigcup_{i' \leq i}\Omega_m^{i'}\setminus \min(\Omega_m^{i})$$ 
is no greater than $$
\left(\min(\Omega_m^{i}) - \frac{\min(\Omega_m^{i})}{(1 + \epsilon)^m} \right)2^m$$
and it follows that 
$$\min(\Omega_m^{i+1}) < \min(\Omega_m^{i}) + \min(\Omega_m^{i})
\left(1 - \frac{1}{(1 + \epsilon)^m} \right)2^m$$
as required.

The general statement \ref{eq:4c} follows by repeated application 
of~\ref{eq:3c}.

To prove \ref{eq:5c} let $n \in \Omega_m^{i}$. Combining
 \ref{eq:2c} and  \ref{eq:4c} yields 
$$ \frac{\Phi(n)}{n} \leq \frac{\min(\Omega_m^{i})}{\min(\Omega_m^{i})} < \left(1 +
  2^m\left(
1 - \frac{1}{(1 + \epsilon)^m}\right)\right)^k(1 + \epsilon)^m$$
establishing the last half of the
inequality. For the first half 
note that $\min(\Omega_m^{i}) \leq
\min(\Omega_m^{i+k})$ 
 and hence, from \ref{eq:2c} and  \ref{eq:4c}
it follows that
$$\frac{n}{\Phi(n)} \leq
\frac{\max(\Omega_m^{i})}{\min(\Omega_m^{i})}
\frac{\min(\Omega_m^{i})}{\min(\Omega_m^{i+k})}
> \frac{1}{(1+ \epsilon)^m} \geq (1 - \epsilon)^m.$$
\end{proof}
\begin{lemma}
  \label{l:d1}
If $p \in \Poset({\mathcal F})$ and $u\in \Naturals$ then there is $f$
such that $(f,m_p,\epsilon_p) \leq p$ and $ u \in \dom(f)$.
\end{lemma}
\begin{proof}
  It suffices to prove this for the case that $u$ is the least 
  integer not in the domain of $f_p$.  Let $i = I_p$ and let $f = f_p
  \cup \Phi\cup \Phi^{-1}$ where $\Phi : \Omega_m^i \to
  \Omega_m^{i+1}$ is an isomorphism. From
  Condition~\ref{eq:A2.1}
of Definition~\ref{d:po} and
  Conclusion~\ref{eq:5c} of  Lemma~\ref{l:cm1}, with $k=1$
it follows that Requirement~\ref{eq:A4} is satisfied. Since
  requirement~\ref{eq:A3} is immediate, this suffices.
\end{proof}
\begin{lemma}
  \label{l:d2}
If $p \in \Poset({\mathcal F})$ and $\epsilon > 0 $ and $m \in \Naturals$ then there is $q \leq p$ such that $\epsilon_q \leq \epsilon$ and $m_q\geq m$.
\end{lemma}
\begin{proof}
  First apply Lemma~\ref{l:d1} to extend 
$p$ to $(f,m_p,\epsilon_p)$ such that 
the domain of $f$ is
  sufficiently large that it is possible to change $m_p$ to $m$ and
  $\epsilon_p$ to $\epsilon$ and still preserve Condition~\ref{eq:A0},
\ref{eq:A0+} and~\ref{eq:A2.1}.
\end{proof}
\begin{lemma}
  \label{l:d3}
If $p \in \Poset({\mathcal F})$ and $\pi\in \Sacks({\mathcal I}_{1/x})$ but
$[\pi]_{{\mathcal I}_{1/x}}$ is not first order definable
from 
 $\{\pi_i\}_{i=0}^\infty/\Fin({\mathcal I}_{1/x})$ 
and $k\in\Naturals$ then there is $q\leq p$ such that
$\sum_{i=k}^\infty\{h(i) :  \pi(f_q(i)) \neq f_q(\pi(i))\} > 1$. 
\end{lemma}
\begin{proof} As a convenience, let $m = m_p$ and $J = I_p$.
\relax From the definition of $J$ it follows that it is possible to choose
 a family of mappings $\{\Phi_{i,j}\}_{i,j\geq J}$ such that
$\Phi_{i,j}:\Omega_m^i\to \Omega_m^j$ is an isomorphism, 
$\Phi_{i,j}\circ\Phi_{j,i}$ is the identity on $\Omega_m^i$
and $\Phi_{i,k}\circ\Phi_{j,i} = \Phi_{j,k}$.
 Also, by appealing to
  Lemma~\ref{l:d1} it may assumed that
if 
$$\delta = \sup_{j \in \Naturals\setminus \dom(f_p)}\left(\sup_{n \leq m_p}
\left | 1 - \frac{\pi_n(j)}{j} \right |\right)$$
then
\begin{equation}
  \label{eq:A2.5}
\left(1 +
  2^m\left(
1 - \frac{1}{(1 + \delta)^m}\right)\right)^6(1 + \delta)^m < 1 + \epsilon_p .
\end{equation}
The following fact will play a role later in the proof but is included
here to explain the significance of the exponent 6 in
 \ref{eq:A2.5} as well as in the indexing to follow.
\begin{claim}\label{f:1}
  Given any $\pi \in Sym(6)$ other than the identity there is $\sigma
  \in Sym(6)$ without fixed points such that $\sigma$ is an
  involution and
  $\sigma$ does not commute with $\pi$.
\end{claim}
Define
 $E_i = \bigcup_{w=0}^5 \Omega_m^{J + 6i+ w}$ for each $ i \in
 \Naturals$. For $\sigma$ an involution in $Sym(6)$ and integers $i$ and $j$
 let $H_{i,j}^\sigma:E_i \to E_{j}$ be the isomorphism
 defined by
$$H_{i,j}^\sigma = \bigcup_{w=0}^5 \Phi_{J + 6i + w, J + 6j +
 \sigma(w)} \cup \Phi_{ J + 6j
 \sigma(w), J + 6i + w} .$$ 
If $\sigma$ is the identity then $H_{i,j}^\sigma$ will be denoted by
$H_{i,j}$ and if $i=j$ then $H_{i,j}^\sigma$ will be denoted by
$H_{i}^\sigma$.

Observe that if $X\subseteq \Naturals$ and $y_i \in E_i$ for $i\in X$ then, 
by \ref{eq:5c} of  Lemma~\ref{l:cm1}, 
\begin{equation}
  \label{eq:Q1}
\sum_{i\in X} 1/y_i < \infty
\text{ if and only if }\sum_{i\in X}\sum_{y\in E_i} 1/y < \infty .  
\end{equation}
This will be used repeatedly in order to restrict the possible
structure of $\pi$.

To begin,
let $W(x)$ be defined so that if $x \in \Omega_m^{i}$ then $\pi(x) \in
\Omega_m^{W(i)}$.
First note that if there exist $x_1$ and $x_2$   in $\Omega_m^{u}$ 
such that $W(x_1) \neq W(x_2)$  then there exists 
 a sequence $(u_1,u_2,
  \ldots , u_k)$ of integers in $m$ such that 
$$x_2 = \pi_{u_1}\circ\pi_{u_2}\circ\ldots
\pi_{u_k}(x_1)$$
and, hence,
 $\pi(\pi_{u_1}\circ\pi_{u_2}\circ\ldots
\pi_{u_k}(x_1)) = \pi(x_2) \neq \pi_{u_1}\circ\pi_{u_2}\circ\ldots
\pi_{u_k}(\pi(x_1))$
because the elements of $\Omega_m$ are closed under the 
permutations $\pi_{u_i}$. 
Therefore, by \ref{eq:Q1}, it may be assumed that
if $Z$ is the set of all $i$ such that there is a pair
$\{x,y\}\subseteq  a \in \Omega_m$ and $a\subseteq E_i$ 
and  $W(x) \neq W(y)$ then $\sum_{i\in Z}\sum_{y\in E_i} 1/y <
\infty$.

Next, let $W'$ be defined such that if $x \in E_i$ then
$\Omega_m^{W(x)}\in E_{W'(x)}$. 
Let $$X = \{i\in\Naturals : (\exists y \in E_i)(\exists x \in E_i)
W'(y) \neq W'(x)\}$$
and suppose that
$\sum_{i\in X}\sum_{y\in E_i} 1/y = \infty$.
Then, using~\ref{eq:Q1},
 it is possible to choose a finite $A\subseteq X$ such that
if $y_a\in E_a$ and $x_a \in E_a$ are chosen for each $a\in A$ such
that $W'(y_a) \neq W'(x_a)$ then $\sum_{a\in A}1/y_a > 1$.
Let $\sigma_a $ be a fixed point free involution of 6 such that 
if $y_a \in \Omega_m^{J + 6i+ w_y}$ and
$x_a \in \Omega_m^{J + 6i + w_x}$ then
$\sigma_a(w_y) = w_x$.
It follows that if $q\leq p$ is the condition defined by
$$f_q = f_p \cup \bigcup_{a\in A}H_a^{\sigma_a}$$
then, by \ref{eq:A2.5} and Conclusion~\ref{eq:5c} of Lemma~\ref{l:cm1}, it follows
that Condition~\ref{eq:A2.1} of Definition~\ref{d:po} is satisfied.
Since each $\sigma_a$ is fixed point free it follows that
Condition~\ref{eq:A3} 
of Definition~\ref{d:po} is also satisfied. 
Therefore $q \leq p$ and it follows from the choice of $A$ that $q$
satisfies the requirements of the lemma.

Hence, 
 it may be assumed that there is a set $Y$ such that
$\bigcup_{i\in \Naturals\setminus Y}E_i \in {\mathcal I}_{1/x}$
and  $W'$ is constant on $E_i$ for each $i\in Y$. 
Let $W''$ be defined on $Y$ such that if $x \in E_i$ then
 $\Omega_m^{W'(x)}\subseteq E_{W''(i)}$. 
Therefore there is a 
partition $Y = Y_0\cup Y_1 \cup Y_2\cup Y_4$ such that 
$W''(Y_i)\cap Y_i = \emptyset$ for each $i \in 3$ and $W''$ is the
identity on $Y_4$. Let $j\in 4$ be
such that $\sum_{i\in Y_j}\sum_{y\in E_i} 1/y = \infty$.

First assume that $j \in 3$.
Let $A\subseteq Y_j$ be such that $\sum_{a\in A}1/y_a > 1$ for any
choice of $y_a\in E_a$ and let $B$ be the image of $A$ under $W''$. 
Choose involutions without fixed points 
$\sigma_a\in Sym(6)$ and $\beta_a\in Sym(6)$ for $a\in A$
such that $$H_{W''(a)}^{\beta_a}\circ\pi\restriction E_a \neq
\pi\circ H_a^{\sigma_a}$$  
and let $q\leq p$ be defined so that
$$f_q = f_p \cup\bigcup_{a\in A}H_a^{\sigma_a} \cup H_{W''(a)}^{\beta_a} $$ 
noting that $q \leq p$ as in the previous case.

Hence, assume that $j = 3$. For $i \in Y_3$ 
let $W_i\in Sym(6)$ be the permutation defined by
$W(x) \in \Omega_m^{W_i(w)}$ if $x \in \Omega_m^{J + 6i + w}$.
 If the set of $i \in Y_3$ such that
$W_i$  is not the identity is not in ${\mathcal I}_{1/x}$ then 
using Claim~\ref{f:1} it is possible to choose
an involution  $\sigma_i\in Sym(6)$ without fixed points such that 
  $\sigma^i$ does not commute with 
$W_i$.
As before,  using \ref{eq:Q1} it is possible to find a finite subset
$A\subseteq Y_3$ such that
setting
$f_q = f_p \cup \bigcup_{a\in A}H_a^{\sigma_a}$ 
suffices.

Therefore, by omitting a set in ${\mathcal I}_{1/x}$, 
it may be assumed that $W_i$ is
the identity for all $i \in Y_3$ and that 
$\bigcup_{i\in \Naturals\setminus Y_3}E_i \in {\mathcal I}_{1/x}$. 
For $\rho \in Sym( \Omega_m^{J})$ and $z\in 6$
let $Y(\rho,z)$ be the set of all $i \in Y_3$
such that
$$\Phi_{J + 6i+ z, J}\circ\pi\circ\Phi_{J, J + 6i
  +z}\restriction 
\Omega_m^{J}  = \rho .$$
 If for each $z \in 6$ there
  is only one $\rho_z \in Sym(\Omega_m^{J})$ such that
$\sum_{i\in Y(\rho_z,z)}\sum_{x\in \Omega_m^{i}}1/x = \infty$ then
$[\pi]_{{\mathcal I}_{1/x}} $ can be defined from
$\{\pi_j\}_{j\in m}$ and $\{\rho_z\}_{z\in 6}$.
 So it may be assumed that it is possible to
choose $z\in 6$ and $\rho_z\in Sym(\Omega_m^{J})$ such that
if $U_0 = Y(\rho_z,z)$ and $U_1 = \Naturals \setminus U_0$
 $$\sum_{i\in U_0}\sum_{x\in \Omega_m^{i}}1/x =\sum_{i\in
   U_1}\sum_{x\in \Omega_m^{i}}1/x= \infty .$$ 
The key point to keep in mind is that if $i \in U_0$ and $j\in U_1$
and $q\leq p$ and $H_{i,j}\subseteq f_q$
 then there is $x\in E_i$ such that $f_q(\pi(x))\neq \pi(f_q(x))$.  
Of course, if $|i-j|$ is too large then there might not be any $q$
such that $ H_{i,j}\subseteq f_q$. The remainder of the argument is
devoted to showing that there are sufficiently many pairs $(i,j)\in
U_0\times U_1$ such that $f_p$ can be extended by $ H_{i,j}$.

To this end, 
let $U = \{n\in U_0 : n+1 \in U_1\}$. For $n \in U$ let $n_0$ be
the  greatest integer  such that the interval
 $[n - n_0 , n]$ is contained in  $U_0$ and 
let $n_1$ be
the  largest integer  such that $[n +1,n+1 + n_1]\subseteq U_1$. Let
$U^*_0$ be the set of all $n\in U$ such that $n_0 \leq n_1$  and
$U^*_1$ be the set of all $n\in U$ such that $n_0 > n_1$ and define
$U'_i = \bigcup_{n\in U^*_i}[n_0,n_1]$ and observe that $Y_3 = U'_0
\cup U'_1$. 
Hence, either  $U'_0$ or $U'_1$ belongs to ${\mathcal I}_{1/x}$ 
  In either case the following argument is similar so
 assume that the former holds.

Let $\bar{\epsilon}$ be so small that
$$\frac{1}{1 - 2^m6\bar{\epsilon}} < 1 +  \epsilon$$
and choose $\delta$ small enough that
$$\frac{\left(1+ 2^mJ((1 +
  \delta)^{6m} - 1)\right)(1 + \delta)^{6m}}{1 - 2^m6(\bar{\epsilon}+ (1 +
  \delta)^{6m} - 1)} .$$ 
Using Lemma~\ref{l:d1}, is may be assumed that
$$1 - \delta <\frac{\pi_i(x)}{x} < 1 + \delta$$
 for all $x
\notin \dom(f_p)$ and $ i \in m$.

There are again two cases to
consider: Either there are infinitely many $n \in U'_0$ such that
$n_0 \leq (1 - \bar{\epsilon})n$ or there are not.
In either case
let $k_n = \max (n-n_0,\lceil n(1 -
\bar{\epsilon})\rceil)$.
Keep in mind that $\min(E_i) = \min(\Omega_m^{J + 6i})$.
Begin by observing, using Conclusion~\ref{eq:2c},
 that if 
$\Omega_m^{j}$ intersects the interval 
$[\lceil\min(E_{k_n})(1 + \delta)^{6m}\rceil,
\min(E_{n})]$ then  $J + 6k_n \leq j \leq J + 6n$.
Hence, 
\begin{equation}
  \label{eq:e19.1}
[\lceil\min(E_{k_n})(1 + \delta)^{6m}\rceil,
\min(E_{n})]\subseteq 
\bigcup_{z = k_n}^n E_z \subseteq \bigcup_{z = \lceil n(1 -
  \bar{\epsilon})\rceil}^n E_z\subseteq \bigcup_{z\in U'_i}
 E_z  
\end{equation}
for $n\in U'_0$.
Furthermore,
\begin{equation}
  \label{eq:e19.2}
 \min(E_{n}) - \lceil\min(E_{k_n})(1 + \delta)^{6m}\rceil
 \leq 
 2^m6(n-k_n)
\end{equation}
It follows that
 $$\min(E_{n}) - \min(E_{k_n}) \leq
 2^m6(n-k_n) +   ( \lceil\min(E_{k_n})(1 + \delta)^{6m}\rceil - 
\min(E_{k_n}))$$
$$ \leq  2^m6(n-k_n) + 2^m(J + 6k_n)((1 + \delta)^{6m} - 1) $$
$$\leq  2^m(6(n - k_n)  + 6k_n((1 + \delta)^{6m} - 1)) + 2^mJ((1 +
\delta)^{6m} - 1)$$  
$$\leq  2^m(6(n - k_n)  + 6n((1 + \delta)^{6m} - 1)) + 2^mJ((1 +
  \delta)^{6m} - 1)$$
$$\leq 2^m6n(\bar{\epsilon}+ (1 + \delta)^{6m} - 1) + 2^mJ((1 +
  \delta)^{6m} - 1) $$ 
and hence, 
$$
\frac{\min(E_n)}{\min(E_{k_n})} - 1
\leq \frac{n}{\min(E_{k_n})}\left(
2^m6(\bar{\epsilon}+ (1 + \delta)^{6m} - 1)
\right) + 2^mJ((1 +
  \delta)^{6m} - 1) .$$
Therefore, using the fact that $n \leq \min(E_n)$,
$$\frac{\min(E_n)}{\min(E_{k_n})}
\leq \frac{1+ 2^mJ((1 +
  \delta)^{6m} - 1)}{1 - 2^m6(\bar{\epsilon}+ (1 + \delta)^{6m} - 1)} \leq
\frac{1 + \epsilon}{(1 + \delta)^{6m}} $$
and so, using Conclusion~\ref{eq:2c}, 
$$\frac{H_{k_n,n}(x)}{x}
\leq \frac{\max(E_n)}{\min{E_{k_n}}} \leq {1 + \epsilon} . $$
Similar reasoning shows that
$$\frac{H_{i, i+k_n}(x)}{x} \leq {1 + \epsilon}  $$ for each $i$ such that $k_n \leq i < n$.
Therefore extending $f_p$ to $f_q$ so that
$f_q = f_p \cup \bigcup_{i= k_n}^{n}H_{i,{i+k_n}}$ 
will satisfy
Conditions~\ref{eq:A3} and~\ref{eq:A4} of Definition~\ref{d:po}.

The only question which remains is whether it is possible to add
enough of these extensions to provide a large witness to
$f_q$ not commuting with $\pi$. 
In the case there is some $K$ such that $k_n = n_0$ for
all  $u\geq K$ it follows that
$\sum_{n \in U'_0\setminus K}\sum_{j=n_0}^n 1/j = \infty$.
Moreover, for each $n \geq K$ and $j$ such that $n_0 \leq j\leq n$
there is some $x_j \in E_{j}$ such that
$\pi(H_{j,{j+n_0}}(x_j)) \neq
H_{j,{j+n_0}}(\pi(x_j))$.
 Hence,
by \ref{eq:Q1}, it follows that 
$\sum_{j \in U'_0\setminus K}1/x_j = \infty$ and so it is possible to
choose
$M$ such that defining $f_q = f_p
=\bigcup_{n=K}^M\bigcup_{i=n_0}^nH_{i,{i+ k_n}}$ works.

In the other case, there is an infinite set $U''\subseteq U'_0$ such
that $k_n = \lceil n(1 -
\bar{\epsilon})\rceil$ for each $n \in U''$. By \ref{eq:e19.1} it
follows that if $n \in U''$ then if  $j \in \lceil u(1 -
\bar{\epsilon})\rceil \leq j \leq u$ then there is some $x_j \in E_{j}$ such that
$\pi(H_{j,{j+n_0}}(x_j)) \neq
H_{j,{j+n_0}}(\pi(x_j))$.
It follows that for $n \in U''$
$$\sum_{i= \lceil n(1-
\bar{\epsilon})\rceil }^n1/x_i \geq 
\sum_{i= \lceil n(1-
\bar{\epsilon})\rceil }^{n}\frac{1}{\max(E_{i})}
\geq\sum_{i= \lceil n(1-
\bar{\epsilon})\rceil }^{n}\frac{1}{\min\left(\Omega_m^{J+ 6i+ 5}\right)(1 + \delta)^{m}}$$
$$\geq\sum_{i= \lceil n(1-
\bar{\epsilon})\rceil }^{n}\frac{1}{ 2^m(J + 6i+ 5)  (1 + \delta)^{m}} $$ 
and elementary calculations show that 
the limit as $n$ increases to infinity of the last term of the inequality is 
$$\frac{1}{2^m6(1 + \delta)^{m}}\ln\left(\frac{1}{1-\bar{\epsilon}}\right) = \gamma >0.$$ 
Now it suffices to choose a finite subset $X
\subseteq U''$ such that $$|X| > 1/\gamma$$ and define 
$f_q
=f_p \cup \bigcup_{n\in X}\bigcup_{i=k_n}^nH_{i,{i+ k_n}}$.
\end{proof}
\begin{theor}\label{t:ca}
  It is consistent that $A({\mathcal I}_{1/x}) = \aleph_1 < 2^{\aleph_0}$.
\end{theor}
\begin{proof}
  Let $V$ be a model where $2^{\aleph_0} > \aleph_1$ and let $V'$ be obtained from $V$ by adding $\aleph_1$ Cohen reals.
Choose permutations $\{\pi_\alpha\}_{\alpha\in\omega_1}$ such that
$\pi_\alpha = \bigcup_{p\in G_\alpha}f_p$ and $G_\alpha\subseteq \Poset(\{\pi_\beta\}_{\beta\in\alpha})$ is generic over
$V[\{\pi_\beta\}_{\beta\in\alpha}]$.
Using Lemmas~\ref{l:d1} and~\ref{l:d2} it follows that $\{\pi_\alpha\}_{\alpha\in\omega_1}$
pairwise almost commute.
Let $G\supseteq \{\pi_\alpha\}_{\alpha\in\omega_1}$ be a maximal
almost Abelian subgroup\footnote{Observe that while the elements of
  the group generated by $\{\pi_\alpha\}_{\alpha\in\omega_1}$ pairwise
  almost commute the same can not be said of the elements of $G$. All
  that can be said of them is that they almost commute modulo the ideal
  ${\mathcal I}_{1/x}$.}
of the subgroup of all 
$\pi\in \Sacks({\mathcal I}_{1/x})/\Fin({\mathcal I}_{1/x})$ which are first
order definable from some finite subset of
$\{\pi_\alpha\}_{\alpha\in\omega_1}$. 
To see that $G$ is maximal in  $\Sacks({\mathcal I}_{1/x})/\Fin({\mathcal
  I}_{1/x})$ suppose that $\pi \in V[\{\pi_\beta\}_{\beta\in\alpha}]$. If  
$\pi$ is first order definable from some finite subset of $\{\pi_\beta\}_{\beta \in\alpha}$ then either
$\pi \in G$ or there is some $\theta \in G$ such that $\{n\in \Naturals : \pi(\theta(n)) \neq \theta(\pi(n))\} \in {\mathcal I}_{1/x}^+$.
On the other hand, if $\pi$ is not first order definable from some finite subset of $\{\pi_\beta\}_{\beta \in\alpha}$ then
by Lemma~\ref{l:d3} and genericity it follows that
$\{n\in \Naturals : \pi(\pi_\alpha(n)) \neq \pi_\alpha(\pi(n))\} \in {\mathcal I}_{1/x}^+$.
\end{proof}

\section{It is possible that ${\mathfrak a}({\mathcal I}) < A({\mathcal
  I})$}
Since it has been shown in Proposition~\ref{p:1} that
$A([\Naturals]^{<\aleph_0}) 
\leq {\mathfrak a}$ it is natural to wonder whether there might not be
a more general result asserting that
 $A({\mathcal I})$ is bounded by $ {\mathfrak
   a}({\mathcal I})$ as defined in Definition~\ref{d:d12}.
It will be shown that no such result holds, at least not in the
generality indicated. 

Fix an increasing sequence of integers ${\mathcal N} =
\{n_i\}_{i=0}^\infty$ such that 
$$\lim_{i\to\infty} \frac{n_{i+1} - n_i}{n_{i+2} - n_{i+1}} = 0$$ and
define $${\mathcal 
  I}({\mathcal N})
= \left\{A\subseteq\Naturals : \lim_{i\to\infty}\frac{|A\cap [n_i,
  n_{i+1})|}{n_{i+1} - n_i} = 0 \right\}.$$
\begin{theor}\label{t:61}
  $A({\mathcal
  I}({\mathcal N})) = 2^{\aleph_0}$.
\end{theor}
\begin{proof}
  To begin the following claim will be established:
  \begin{claim}\label{c:11}
    If $g\in \Sacks({\mathcal
  I}({\mathcal N}))$ then there is $B\in{\mathcal
  I}({\mathcal N})$ such that if $j\in [n_i, n_{i+1})\setminus B$ then 
$g(j)\in [n_i, n_{i+1})$.
  \end{claim}
  \begin{proof}
    Let $B^+ = \bigcup_{i=0}^\infty\{n_i \leq n <  n_{i+1} :  g(n)
    \geq n_{i+1}\}$ and 
    let $B^- = \bigcup_{i=0}^\infty\{n_i \leq n <  n_{i+1} :  g(n)
    < n_{i}\}$. If $B^+\cup B^- \in {\mathcal
  I}({\mathcal N})$ then the claim is proved.
To begin suppose $B^+\notin {\mathcal
  I}({\mathcal N})$. Choose $\epsilon > 0$ and an  infinite $Y\subseteq
    \Naturals$ such that $$\frac{|B^+\cap [n_i, n_{i+1})|}{n_{i+1} -
    n_i} \geq \epsilon$$ for each $i \in Y$. By thinning out $Y$ it
    may also be assumed that if $i$ and $j$ belong to $Y$ and $i < j$
    then $g(B^+\cap [n_i,n_{i+1})\subseteq n_j$. It follows that $g(B^+)\cap
    [n_{i+1}, n_j) = g(B^+\cap [n_i,n_{i+1})\cap
    [n_{i+1}, n_j)$. Therefore, if $i < k < j$ then 
$$\frac{|g(B^+)\cap [n_k, n_{k+1})|}{n_{k+1} -
    n_k} \leq \frac{n_{i+1}) - n_i}{n_{k+1} -
    n_k}$$
and so $g(B^+) \in {\mathcal
  I}({\mathcal N})$ contradicting that $g \in \Sacks({\mathcal
  I}({\mathcal N}))$. A similar argument applied to $g^{-1}$ deals
with $B^-$.
  \end{proof}
Given a subgroup $G\subseteq \Sacks({\mathcal
  I}({\mathcal N}))/\Fin({\mathcal
  I}({\mathcal N}))$ let $g \in G$ be different from the identity. Hence
$\{n :f(n) \neq n\} \in {\mathcal
  I}({\mathcal N})^+$ and, using the Claim, there is $B \in {\mathcal
  I}({\mathcal N})$ such that
 if $j\in [n_i, n_{i+1})\setminus B$ then 
$g(j)\in [n_i, n_{i+1})$. Let $g'$ be a permutation such that
$g'\restriction \Naturals \setminus B = g\restriction \Naturals \setminus B$
and $g'\restriction [n_i, n_{i+1})$ is a permutation of $[n_i,
  n_{i+1})$ for each $i$. (This is possible since $g\restriction
   [n_i, n_{i+1})\setminus B \to [n_i, n_{i+1})$ is one-to-one.)
Note that $g' \in \Sacks({\mathcal
  I}({\mathcal N}))/\Fin({\mathcal
  I}({\mathcal N}))$. For any $Z\subseteq \Naturals$ let $g_Z$ be defined
  by 
 $$g_Z(j) = 
 \begin{cases}
   g'(j) & \IF j \in [n_i, n_{i+1}) \AND i\in Z\\
j & \IF j \in [n_i, n_{i+1}) \AND i\notin Z
 \end{cases}
$$ and note that $g_Z \in \Sacks({\mathcal
  I}({\mathcal N}))/\Fin({\mathcal
  I}({\mathcal N}))$ for each $Z$. 

Now, suppose that  $h \in \Sacks({\mathcal
  I}({\mathcal N}))/\Fin({\mathcal
  I}({\mathcal N}))$ and use the claim to find $C\in {\mathcal
  I}({\mathcal N})$ such that
 $j\in [n_i, n_{i+1})\setminus B$ then 
$g(j)\in [n_i, n_{i+1})$. Let $D$ be such that $h(g(j))= g(h(j))$ for
  $j \in \Naturals \setminus D$.
Then let $E = C\cup h(B) \cup h^{-1}(B)\cup D$
  and note that $E \in {\mathcal
  I}({\mathcal N})$ since $h \in  \Sacks({\mathcal
  I}({\mathcal N}))/\Fin({\mathcal
  I}({\mathcal N}))$. Now observe that
$g_Z(h(j)) = h(g_Z(j))$ for each $j\in \Naturals \setminus E$.
To see this let $j \in [n_i, n_{i+1})$ and suppose first that $ i\in
  Z$.
In this case $g_Z(j) = g'(j) = g(j)$. Furthermore, 
since $j\notin C$, $h(j) \in [n_i, n_{i+1})$ and hence
 $g(h(j)) = g_Z(h(j)$. Since $j\notin D$ it follows that
$h(g(j)) = g(h(j)$ and, hence,
In this case $h(g_Z(j)) = g_Z(h(j))$.
 If $ i\notin
  Z$ then $g_Z(j) = j$ and, since $j\notin D$, $h(j) \in [n_i, n_{i+1})$.
Therefore, $g_Z(h(j)) = h(j)$ and so  $g_Z(h(j)) = h(g_Z(j))$.
\end{proof}

\begin{theor}
  $\mathfrak a({\mathcal
  I}({\mathcal N})) \leq \mathfrak a$.
\end{theor}
\begin{proof}
  Let ${\mathcal A}$ be a maximal almost disjoint family of size
  $\mathfrak a$. For $A\in {\mathcal A}$ define $A^* = \bigcup_{i\in
  A}[n_i, n_{i+1})$ and let ${\mathcal A}^*= \{A^* : A\in {\mathcal
  A}\}$.
Then ${\mathcal A}^*$ is maximal in ${\mathcal P}(\Naturals)/{\mathcal
  I}({\mathcal N})$.
\end{proof}

\section{Questions}
\begin{quest}
  Can the lower bound $A(\Sacks/\Fin) \geq \mathfrak p$ of
Theorem~\ref{t:pt} be improved?
\end{quest}
For any function $h : \Naturals \to \Reals$ one can define
the {\em summable}  ideal ${\mathcal I}_h$  to be the set of all
$X\subseteq \Naturals$ 
such that $\sum_{x\in X}h(x) < \infty$.  
Observe that it is possible to modify the proof of Theorem~\ref{t:61}
in order to replace the ideal ${\mathcal
  I}({\mathcal N})$ by a summable ideal. In particular, let
$\{n_i\}_{i=0}^\infty$ be an  increasing sequence of integers defined by
$n_{i+1} - n_i = n_i^3$ and let $h$ be defined by
$h(j) = n_i^{-3}$ if $n_i \leq j < n_{i+1}$. If $g\in \Sacks({\mathcal
  I}_h)$ and $B^+$ and $B^-$ are defined  as in the proof of
Theorem~\ref{t:61} then it is easy to see that $\sum_{j\in B^+\cap
  n_i} h(g(j)) \leq |B^+\cap
  n_i|n_i^{-3} \leq n_i^{-2}$ and hence Claim~\ref{c:11} still holds
  as does the remainder of the argument of Theorem~\ref{t:61}.
Hence $A({\mathcal I}_h)= 2^{\aleph_0}$. This motivates the following
question.

\begin{quest}
  For which functions $h$ is it possible to improve Theorem~\ref{t:ca}
to show that $A({\mathcal I}_h)= \aleph_1 < 2^{\aleph_0}$ in the model
obtained by adding $\aleph_1$ Cohen reals? 
\end{quest}
\begin{quest}
  Are there functions $h$ and $g$ such that it is consistent that
$A({\mathcal I}_h) <A({\mathcal I}_g) < 2^{\aleph_0}$?
\end{quest}
\begin{quest}
  Is it possible to characterize the summable ideals ${\mathcal I}_h$ such that
$A({\mathcal I}_h)= 2^{\aleph_0}$? Can the same be done for the $F_\sigma$
ideals or all analytic ideals? 
\end{quest}
\begin{quest}
  Can Theorem~\ref{t:ca} be improved to show that it is consistent
  with set theory that $2^{\aleph_0} > \aleph_1$ yet there is an
  almost commuting subgroup of $\Sacks$ of cardinality $\aleph_1$
 which is maximal with respect
  to commuting modulo ${\mathcal I}_{1/x}$ ? Does this hold in the
  Cohen model of Theorem~\ref{t:ca}?
\end{quest}
The methods of Sections~4 and~5 require that the subgroups constructed
contain many involutions. While the methods can be modified to produce groups
with no elements of order $k$ for a fixed $k$, the following questions
seem more subtle.
\begin{quest}
  Can Theorem~\ref{tt:ml} be modified to assert that
it is consistent that there is a maximal, Abelian, torsion free
subgroup of $\Sacks/\Fin$ of
size $ \aleph_1$ and  $\aleph_1< {\mathfrak a}$ ? 
\end{quest}
\begin{quest}
  Can Theorem~\ref{t:ca} be modified to assert that
it is consistent that there is a maximal, Abelian, torsion free
subgroup of $\Sacks({\mathcal I}_{1/x})/\Fin({\mathcal I}_{1/x})$ of
size
$ \aleph_1$ and  $\aleph_1< 2^{\aleph_0}$ ? 
\end{quest}


\end{document}